\documentclass{amsart}
\usepackage{amssymb, amssymb, amsmath, amsthm, amscd, latexsym, mathrsfs, tgcursor}
\input xy
\xyoption{all}
\usepackage{times}
\usepackage{enumerate}

\newtheorem{thm}{Theorem}[section]
\newtheorem{cor}[thm]{Corollary}
\newtheorem{lem}[thm]{Lemma}

\newtheorem{prop}[thm]{Proposition}
\newtheorem{hypa}{Hypothesis}
\renewcommand{\thehypa}{A\,\,}
\newtheorem{hypc}{Hypothesis}
\renewcommand{\thehypc}{C\,\,}
\newtheorem{hypb}{Hypothesis}
\renewcommand{\thehypb}{B\,\,}
\theoremstyle{definition}

\newtheorem{rem}[thm]{Remark}
\newtheorem{exa}[thm]{Example}

\numberwithin{equation}{section}

\newcommand{\lra}{\longrightarrow}
\newcommand{\co}{\colon\!}

\newcommand{\smin}{\smallsetminus}
\newcommand{\pt}{{*\,}}

\newcommand{\id}{\textup{id}}

\newcommand{\holim}{\textup{holim}}
\newcommand{\hocolim}{\textup{hocolim}}
\newcommand{\colim}{\textup{colim}}
\newcommand{\hofiber}{\textup{hofiber}}
\newcommand{\mor}{\textup{mor}}
\newcommand{\map}{\textup{map}}
\newcommand{\nat}{\textup{nat}}
\newcommand{\Rnat}{\textup{Rnat}}
\newcommand{\Rmap}{\textup{Rmap}}
\newcommand{\reg}{\mathscr R}
\newcommand{\ho}{\textup{h}}
\newcommand{\simp}{\textup{simp}}
\newcommand{\op}{\textup{op}}

\newcommand{\TOP}{\textup{TOP}}
\newcommand{\STOP}{\textup{STOP}}
\newcommand{\GL}{\textup{GL}}
\newcommand{\G}{\textup{G}}
\newcommand{\SOr}{\textup{SO}}
\newcommand{\Or}{\textup{O}}
\newcommand{\aut}{\textup{Aut}}
\newcommand{\diff}{\textup{diff}}
\newcommand{\topo}{\textup{top}} 
\newcommand{\Ofun}{\textup{\textsf{Bo}}}  
\newcommand{\Tfun}{\textup{\textsf{Bt}}}
\newcommand{\TOfun}{\textup{\textsf{Bto}}}
\newcommand{\Gfun}{\textup{\textsf{Bg}}}
\newcommand{\Sfun}{\textup{\textsf{K}}}
\newcommand{\Afun}{\textup{\textsf{A}}}
\newcommand{\Efun}{\textup{\textsf{E}}}
\newcommand{\Ffun}{\textup{\textsf{F}}}
\newcommand{\Hfun}{\textup{\textsf{H}}}
\newcommand{\Cfun}{\textup{\textsf{C}}}
\newcommand{\Dfun}{\textup{\textsf{D}}}
\newcommand{\Mfun}{\textup{\textsf{M}}}
\newcommand{\Xfun}{\textup{\textsf{X}}}

\newcommand{\sY}{\mathscr Y}
\newcommand{\sJ}{\mathscr J}
\newcommand{\sE}{\mathscr E}

\newcommand{\NN}{\mathbb N}
\newcommand{\RR}{\mathbb R}
\newcommand{\ZZ}{\mathbb Z}
\newcommand{\QQ}{\mathbb Q}

\newcommand{\holimsub}[1]{\underset{#1}{\holim} \,}
\newcommand{\hocolimsub}[1]{\underset{#1}{\hocolim} \,}
\newcommand{\colimsub}[1]{\underset{#1}{\colim} \,}

\begin{document}




\title{Rational Pontryagin classes and functor calculus}

\author{Rui Reis
\and
Michael Weiss}

\address{Mathematisches Institut, Universit\"{a}t M\"{u}nster, Einsteinstr. 62, 48149 M\"{u}nster, Germany.}
\email{m.weiss@uni-muenster.de}


\maketitle


\begin{abstract}
It is known that in the integral cohomology of $BS\Or(2m)$, the square of the Euler class is the same as the Pontryagin class
in degree $4m$. Given that the Pontryagin classes extend \emph{rationally} to the cohomology of $B\STOP(2m)$, it is reasonable to ask whether
the same relation between the Euler class and the Pontryagin class in degree $4m$ is still valid in the rational cohomology of
$B\STOP(2m)$. In this paper we reformulate the hypothesis as a statement in differential topology,
and also in a functor calculus setting.
\end{abstract}

\section{Introduction}
Let $\Or(n)$ be the orthogonal group of $\RR^n$ and $\TOP(n)$ the group of homeomorphisms from $\RR^n$ to itself,
viewed as a topological or simplicial group. The rational cohomology of $B\Or(n)$ is a polynomial ring generated by the
Pontryagin classes $p_1,p_2,\dots,p_{\lfloor n/2\rfloor}$ where $|p_i|=4i$. The rational cohomology of $B\SOr(n)$ has the
same description for odd $n$, while for even $n$ it can be described as
\[ \QQ[p_1,p_2,\dots,p_{n/2},e]\big/(p_{n/2}-e^2) \]
where $e$ is the Euler class in degree $n$. By contrast the rational cohomology rings of $B\TOP(n)$ and the orientable variant
$B\STOP(n)$ are not well understood. Using the work of Waldhausen on $h$-cobordisms and algebraic $K$-theory,
Farrell and Hsiang \cite{FarrellHsiang} calculated $H^*(B\TOP(n);\QQ)$ in the range
$*<4n/3$ approximately, deducing that for odd $n\gg 0$ the
inclusion $B\Or(n)\to B\TOP(n)$ does not induce an isomorphism in rational
cohomology. Watanabe and Sakai-Watanabe~\cite{Watanabe1, Watanabe2, SakaiWatanabe} also show that the
rational cohomology of $B\TOP(n)$ is in many cases much larger that the rational cohomology of $B\Or(n)$. These
calculations are typically formulated as results about spaces of smooth structures on $n$-disks or $n$-spheres.
To make the connection with $B\TOP(n)$, use smoothing theory as elaborated further on in this introduction and in
section~\ref{sec-orcalsmooth}.

\smallskip
The stable case is completely understood. Indeed the inclusion of $B\Or$ into $B\TOP$ is a rational homotopy
equivalence. This is based on celebrated work of Thom, Novikov, Kirby-Siebenmann and others. It follows that
$H^*(B\TOP;\QQ)$ is a polynomial ring with one generator $p_i$ in each degree $4i$, where $i>0$. By restriction,
the $p_i$ are also defined in $H^*(B\TOP(n);\QQ)$ for every (finite) $n$. Less dramatically, the Euler class
$e$ in $H^n(B\STOP(n);\QQ)$ can be defined for even $n$ using the spherical fibration associated to the universal $\STOP(n)$-bundle.
We can therefore ask whether the relation $e^2=p_{n/2}$ holds in $H^*(B\STOP(n);\QQ)$ for even $n$ as it does in $H^*(B\SOr(n);\QQ)$.

\begin{hypa} \label{hypa} For all positive even integers, the equation $e^2=p_{n/2}$ holds in the cohomology ring
$H^{2n}(B\STOP(n);\QQ)$.
\end{hypa}

Hypothesis \thehypa for a specific $n$ implies easily that $p_{n/2}=0 \in H^{2n}(B\TOP(m);\QQ)$ for $m<n$.
See remark~\ref{rem-eutopo}.

\medskip
We are cautious with regard to hypothesis \thehypa by calling it a hypothesis rather than a conjecture.
The hypothesis comes from functor calculus, specifically the orthogonal calculus which has claims to
be a homotopical theory of characteristic classes. In orthogonal calculus, we consider continuous functors from a
certain category $\sJ$ to spaces. The objects of $\sJ$ are finite dimensional real vector spaces
with inner product and the morphisms are linear (injective) maps respecting the inner product.
Orthogonal calculus applied to the functor $V\mapsto B\Or(V)$,
where $V$ is a finite dimensional real vector space with inner product, does a good job in re-constructing the Euler classes
and Pontryagin classes of vector bundles as structural features of the Taylor tower of the functor. The question arose whether the
orthogonal calculus applied to the functor $V\mapsto B\TOP(V)$ can do an equally good job in reconstructing the
rational Pontryagin classes for bundles with structure group $\TOP(n)$. Hypothesis \thehypa would imply that the answer is yes.
More precisely we show that hypothesis \thehypa is equivalent to the following. Let $\Ofun$ and $\Tfun$ be the functors given by
$V\mapsto B\Or(V)$ and $V\mapsto B\TOP(V)$, respectively.

\begin{hypc}  The inclusion
$\Ofun\to \Tfun$ admits a rational left inverse (up to weak equivalence).
\end{hypc}

Using the Taylor towers of $\Ofun$ and $\Tfun$ we are able to reformulate hypothesis \thehypc as a statement about a map
of spectra with action of $\Or(2)$. This allows us to draw some surprising conclusions regarding hypothesis \thehypa. It turns out,
for example, that if hypothesis \thehypa is wrong for some even $n$, then it is wrong for almost all even $n$.

\medskip
The groups $\TOP(n)$ as structure groups of fiber bundles with fiber $\RR^n$ are important in differential topology
because they are ingredients in smoothing theory. Among several geometric hypotheses equivalent to hypothesis \thehypa,
we find the following particularly convenient for comparison. Let $\mathscr R=\mathscr R(n,2)$ be the space of
smooth regular (=nonsingular) maps from $D^n\times D^2$ to $D^2$ which agree with the projection $D^n\times D^2\to D^2$
on and near the boundary. For $f\in \mathscr R$, the derivative $df$ can be viewed as a map from $D^n\times D^2$ to the based space
$Y$ of surjective linear maps from $\RR^{n+2}$ to $\RR^2$, taking all of $\partial(D^n\times D^2)$ to the base point.
(The space $Y$ can be identified with the space of \emph{injective} linear maps $\RR^2\to \RR^{n+2}$ by taking transposes.
It can also be described as the coset space $\GL_{n+2}(\RR)/\GL_n(\RR)$ and it is homotopy equivalent to
the Stiefel manifold $\Or(n+2)/\Or(n)$ of orthonormal $2$-frames in $\RR^{n+2}$.)
So we have a map
\[  \nabla\co \mathscr R \lra \Omega^{n+2}Y~;~\nabla(f)=df~.\]
If $n$ is even and $\ge 4$, then the base point component of the
target of this map is rationally an Eilenberg-MacLane space $K(\QQ,n-3)$, but there is a nontrivial finite $\pi_0$\,.
The map $\nabla$ is a based $S^1$-map,
with $S^1\cong \SOr(2)$ acting by a form of conjugation on source and target, fixing the base point in each case.
To suppress the nontrivial $\pi_0$ in the target, compose with the inclusion $\Omega^{n+2}Y\to \Omega^{n+2}(Y_\QQ)$,
where $Y_\QQ$ is the rationalization of the path-connected space $Y$. Now
$\Omega^{n+2}(Y_\QQ)$ is path connected, and it is an honest $K(\QQ,n-3)$. We still write
\[ \nabla\co \mathscr R \lra \Omega^{n+2}(Y_\QQ) \]
for the composition of the above $\nabla\co  \mathscr R \lra \Omega^{n+2}Y$ with $\Omega^{n+2}Y\to \Omega^{n+2}(Y_\QQ)$.
This allows us to write
\[  [\nabla]\in H^{n-3}_{S^1}(\mathscr R(n,2),\pt;\QQ) \]
using Borel cohomology.

\begin{hypb} \label{hypb} For all even $n\ge 4$, the class $[\nabla]\in H^{n-3}_{S^1}(\mathscr R(n,2),\pt;\QQ)$ is zero.
\end{hypb}

\smallskip
\emph{Acknowledgments.} We wish to thank Ilan Barnea, David Barnes and Pedro Boavida Brito
for advice related to model categories.
Due partly to their influence, we have gradually adopted more and more model category language in writing this paper.
A sense of guilt remains that we did not go far enough.

\smallskip
\emph{Organisation.} In the short section~\ref{sec-easy} we show that \thehypc implies \thehypa for all $n\ne 4$, and
that \thehypa for a specific $n$ implies \thehypb for the same $n$. We set up functor calculus
machinery in sections~\ref{sec-ortho}
and~\ref{sec-homoext}. In section~\ref{sec-split} we introduce a spectrum variant of \thehypc and show that
this is equivalent to \thehypc.
We also pick up the remaining case $n=4$ in \thehypc$\Rightarrow$\thehypa.
Then we show in section~\ref{sec-orcalsmooth} that \thehypb for infinitely many even $n$ implies the
spectrum variant of \thehypc. In section~\ref{sec-pessimist}, as an afterthought and a concession to
pessimism, we formulate weaker versions of hypotheses \thehypa, \thehypb and \thehypc
and prove their equivalence.

\section{Easy implications} \label{sec-easy}
We start by stating hypothesis \thehypc and related definitions in detail.
The objects of $\sJ$ are finite dimensional real
vector spaces $V$ with a (positive definite)
inner product. A morphism in $\sJ$ from $V$ to $W$ is a linear map respecting the inner product.
It is therefore necessarily injective. We regard $\mor_\sJ(V,W)$ as a space, so that
$\sJ$ is enriched over spaces. As a rule, not always rigorously kept, we mean by a \emph{space}
an object of $\bf T$, the category of compactly generated weak Hausdorff spaces. Similarly, a based space
is an object of ${\bf T}_\pt$. See \cite[Def. 2.4.21, Cor. 2.4.26]{Hovey}. In particular, continuous
functors from $\sJ$ to based spaces are enriched functors from $\sJ$
to ${\bf T}_\pt$~.

\smallskip
We need to discuss the meaning of \emph{rational weak homotopy equivalence}
before we use the concept. It is agreed that a map between connected nilpotent
spaces \cite{HilMisRo} is a rational weak homotopy equivalence if and only if it induces an isomorphism
in rational singular cohomology. If $f\co X\to Y$ is a based map between based connected nilpotent spaces
which is a rational weak homotopy equivalence, then $\Omega f\co \Omega X\to \Omega Y$
restricted to base point components is a rational homotopy equivalence.
Unfortunately, some spaces which are important to us here are not nilpotent.
The space $B\Or(n)$ for even $n>0$ is not nilpotent because $\pi_1B\Or(n)$
acts very nontrivially on $\pi_nB\Or(n)$\,. \newline
We are led to the following compromise. Let
$K$ be a discrete group. Let $f\co X\to Y$ be a map between based path-connected spaces
\emph{over the space $BK$}. That is, $X$ and $Y$ are equipped with reference maps to the
classifying space $BK$ and $f$ respects the reference maps. Suppose that
the path components of $\hofiber[X\to BK]$ and
$\hofiber[Y\to BK]$ are nilpotent spaces. We shall say that $f$ is a rational
weak homotopy equivalence \emph{over $BK$} if the induced map
from $\hofiber[X\to BK]$ to $\hofiber[Y\to BK]$
induces a bijection on $\pi_0$ and, on path components,
is a rational weak equivalence between nilpotent spaces. Note that these homotopy
fibers can also be described as covering spaces of $X$ and $Y$, respectively. \newline
In the cases where $K$ is the trivial group, we shall continue to use the expression
\emph{rational weak homotopy equivalence}.
But very often, in the examples that we shall encounter, the group $K$ is $\ZZ/2$.
In particular, the rational left inverse in hypothesis \thehypc is a natural transformation
$\Tfun\to\Efun$ of continuous functors from $\sJ$ to ${\bf T}_\pt$ such that the composition
$\Ofun\to\Tfun\to\Efun$ specializes to a rational weak equivalence
$\Ofun(V)\to \Tfun(V)\to \Efun(V)$ over $B\ZZ/2$ for every $V$ in $\sJ$.
(We are therefore assuming that $\Efun$ comes equipped with a natural transformation to
the constant functor with value $B\ZZ/2$.)

\begin{rem} In many places in this paper, and especially in this section,
we use the following principle. Let $X$ be a based space, let $A$ be an abelian group and
$s$ a positive integer. Then a singular cohomology class $z$ in $H^r(X,*;A)$ determines a
singular cohomology class $\Omega^sz$ in $H^{r-s}(\Omega^sX,*;A)$. There are many difficult constructions for this,
but the following is easy. Assume without loss of generality
that $X$ has the homotopy type of a based CW-space. The class $z$ is represented by a based map from $X$
to an Eilenberg-MacLane space; by applying $\Omega^s$ to that map, we obtain a map from
$\Omega^sX$ to an Eilenberg-MacLane space. \qed
\end{rem}

\begin{prop} \label{prop-CtoA} Hypothesis \emph{\thehypc} implies the cases $n\ne 4$ of hypothesis \thehypa.
\end{prop}

\proof There is a weakening $\textup{\thehypc\!\!}^\delta$ of \thehypc
where we restrict from $\sJ$ to the subcategory $\sJ^\delta$ which has objects $\RR^0,\RR^1,\RR^2,\dots$
and, as morphisms, only the standard inclusion maps $\RR^i\to \RR^j$ for $i\le j$. We shall prove that $\textup{\thehypc\!\!}^\delta$
implies the cases $n\ne 4$
of hypothesis \thehypa.
Suppose therefore that we have a natural transformation
\[ \Tfun \to \Efun \]
of functors on $\sJ^\delta$ such that the composition $\Ofun(\RR^n)\to \Tfun(\RR^n)\to \Efun(\RR^n)$
is a rational homotopy equivalence over $B\ZZ/2$ for every $n$.
\newline
Let $n$ be even, $n\ne 4$. We can pass to double covers
\[ B\SOr(n)\to B\STOP(n) \to \Efun^\natural(\RR^n)~. \]
Therefore we can talk about induced homomorphisms
\[
\xymatrix{ H^n(B\SOr(n);\QQ) & \ar[l] H^n(B\STOP(n);\QQ) & \ar[l] H^n(\Efun^\natural(\RR^n);\QQ). }
\]
These homomorphisms
are isomorphisms because $B\Or(n)\to B\TOP(n)$ is rationally $(n+2)$-connected.
See remark~\ref{rem-hiconn} below. It follows that the element of 
\[ H^n(\Efun^\natural(\RR^n);\QQ) \]
which corresponds to the Euler class
in $H^n(B\SOr(n);\QQ)$ maps to the Euler class in \linebreak $H^n(B\STOP(n);\QQ)$. Also, the element of
$H^{2n}(\Efun^\natural(\RR^n);\QQ)$ which corresponds to the Pontryagin class
in $H^{2n}(B\SOr(n);\QQ)$ maps to the Pontryagin class in $H^{2n}(B\STOP(n);\QQ)$. This follows from the fact that
the homomorphisms
\[
\xymatrix{ H^{2n}(B\Or;\QQ) & \ar[l] H^{2n}(B\STOP;\QQ) & \ar[l] H^{2n}(\Efun^\natural(\RR^\infty);\QQ) }
\]
are isomorphisms, with $\Efun^\natural(\RR^\infty)=\hocolim_n \Efun(\RR^n)$. Therefore the equation
\[ e^2=p_{n/2} \]
in $H^{2n}(B\SOr(n);\QQ)\cong H^{2n}(\Efun^\natural(\RR^n);\QQ) $ implies the same in $H^{2n}(B\STOP;\QQ)$
by pullback from $\Efun^\natural(\RR^n)$. \qed

\begin{rem} \label{rem-hiconn} It is stated in \cite[Essay V, 5.0.(4)]{KirbySiebenmann} that
$\pi_k(\TOP/\Or,\TOP(n)/\Or(n))$ is zero
if $k\le n+2$ and $n\ge 5$. Since $\TOP/\Or$ is rationally contractible, this implies that 
\[ \pi_{k-1}(\TOP(n)/\Or(n))\otimes\QQ=0 \]
under the same assumptions; equivalently, $\pi_k(B\TOP(n),B\Or(n))\otimes\QQ=0$.
When $n=2$ the inclusion $B\SOr(n)\to B\STOP(n)$ is a homotopy equivalence.
See \cite[Essay V, 5.0.(7)]{KirbySiebenmann}.
\end{rem}

\begin{prop} \label{prop-AtoB} Hypothesis~\ref{hypa} for a specific even integer $n\ge 4$
implies hypo\-thesis~\ref{hypb} for the same $n$.
\end{prop}

\proof We begin with a few words on hypothesis~\ref{hypb}.
The space $Y$ of surjective linear maps from $\RR^n\times\RR^2$ to $\RR^2$ fits into a homotopy fiber sequence
\[   S^n\to Y\to S^{n+1}~. \]
Up to homotopy equivalences this is the unit tangent bundle of $S^{n+1}$. If $n$ is even, then the projection 
$Y\to S^{n+1}$ has a section, so $\pi_0(\Omega^{n+2}Y)\cong\pi_{n+2}(Y)\cong\ZZ/2\oplus\ZZ/2$. Also, it
follows easily that the base point component of $\Omega^{n+2}Y$ is rationally homotopy equivalent to an Eilenberg-MacLane
space $K(\QQ,n-3)$. (We are assuming that $n$ is even and $n\ge 4$.)
The fact that $\Omega^{n+2}Y$ is not path connected is irritating.
To remedy this the recommendation was to (re-)define $\nabla$ in hypothesis \thehypb
as the composition of $\mathscr R \lra \Omega^{n+2}Y~;~f\mapsto df$
with the inclusion $\Omega^{n+2}Y\to \Omega^{n+2}(Y_\QQ)$. Here $\Omega^{n+2}(Y_\QQ)$
is path connected. (The idea of simply replacing $\Omega^{n+2}Y$
by its base point component is not a good alternative because we have not shown that
the map $\mathscr R \lra \Omega^{n+2}Y~;~f\mapsto df$ lands in the base point
component of $\Omega^{n+2}Y$.) \newline
We work in a certain model category of based $S^1$-spaces. The based spaces are objects of ${\bf T}_\pt$\,. A based $S^1$-map
$X\to Y$ is a weak equivalence if the underlying map ($S^1$-actions suppressed) is a weak homotopy equivalence.
It is a fibration if the underlying map ($S^1$-actions suppressed) is a Serre fibration. The model category
determines a so-called \emph{homotopy category}, by a procedure which amounts to formal inversion of the
weak equivalences but is more explicit than that; see \cite[Def.~2.13]{GoerssSch}.  \newline
Let $Z=\TOP(n+2)/\TOP(n)$, a variant of $Y=\GL_{n+2}(\RR)/\GL_n(\RR)$.
The idea is to show that if hypothesis \thehypa holds for our specific $n$, then the based map
\[ \Omega^{n+2}Y_\QQ \lra \Omega^{n+2}Z_\QQ \]
admits a left inverse in the homotopy category of based $S^1$-spaces.
To construct such a left inverse we ought to construct a class
\begin{equation}\label{eqn-classz}  z\in H^{n-3}_{S^1} (\Omega^{n+2}Z_\QQ,\pt;\QQ) \end{equation}
which maps to a generator of $H^{n-3}_{S^1}(\Omega^{n+2}Y_\QQ,\pt;\QQ)\cong\QQ$.
In fact, in order to achieve $S^1$-invariance, we shall construct (after the ongoing preview) a class
\begin{equation} \label{eqn-classz'}
 z'\in H^{n-1}(\Omega^n(\TOP/\TOP(n));\QQ) \end{equation}
and note that the appropriate action of $S^1$ on $\Omega^n(\TOP/\TOP(n))$ is trivial in a weak sense;
i.e., the based $S^1$-space $\Omega^n(\TOP/\TOP(n))$ is weakly equivalent to a based $S^1$-space with
a trivial action of $S^1$.
(Think of $\Omega^n(\TOP/\TOP(n))$ as $\Omega^n(\TOP(n+2+N)/\TOP(n))$
for very large $N$. The group $S^1$ acts via conjugation by rotations
of the summand $\RR^2$ in $\RR^n\oplus\RR^2\oplus\RR^N$. It acts trivially on the $n$-fold loop coordinates.)
Therefore $z'$ can be regarded as a class in
\[ H^{n-1}_{S^1}(\Omega^n(\TOP/\TOP(n)),\pt;\QQ). \]
We apply $\Omega^2$ to obtain
\[ z''\in H^{n-3}_{S^1}(\Omega^{n+2}(\TOP/\TOP(n)),\pt;\QQ) \]
and we go from there to $H^{n-3}_{S^1} (\Omega^{n+2}Z,\pt;\QQ)$ by restriction,
and project to the appropriate direct summand to get $z$ as in~(\ref{eqn-classz}).
Here the action of $S^1$ on the two loop
coordinates is nontrivial, but that does not stop us. \newline
Continuing with the preview: now that we have a homotopy left inverse
for the $S^1$-map \linebreak $\Omega^{n+2}Y_\QQ \lra \Omega^{n+2}Z_\QQ$\,, it only remains to show that the composition
\begin{equation} \label{eqn-lose}
\xymatrix{
\reg(n,2) \ar[r]^-{\nabla} & \Omega^{n+2}Y \ar[r] & \Omega^{n+2}Z
}
\end{equation}
represents the zero morphism in the homotopy
category of based $S^1$-spaces. This can be seen as follows. We place
$\reg(n,2)$ into a homotopy fiber sequence
\[ \reg(n,2)\lra  \Omega^2B_1 \lra B_2~. \]
Here $B_1$ is a classifying space for bundles where the fibers are smooth compact contractible $n$-manifolds
whose boundary is identified with the boundary of $D^n$. Similarly $B_2$ is a classifying space for bundles where the fibers are
smooth compact contractible $(n+2)$-manifolds (with corners) whose boundary is identified with the boundary of
$D^n\times D^2$. There are topological variants of $\reg(n,2)$ and of $B_1$ and $B_2$, related by a similar
homotopy fiber sequence. For example the topological variant
of $B_1$ is a classifying space for bundles where the fibers are compact contractible $n$-manifolds
whose boundary is identified with the boundary of $D^n$. The Alexander trick and the affirmed Poincar\'e conjecture
together prove contractibility of the topological variants of $B_1$ and $B_2$; hence the topological variant
of $\reg(n,2)$ is also contractible. The composition in~(\ref{eqn-lose}) factors through the topological variant of $\reg(n,2)$
and is therefore zero in the homotopy
category of based $S^1$-spaces.
\newline
End of preview --- now for the construction of the class $z'$ in~(\ref{eqn-classz'}).
We make a chase through the diagram of (rational) cohomology groups
\begin{equation} \label{eqn-AtoB} \begin{split}
\xymatrix@R=15pt{ H^{n}(\Omega^nB\TOP,\Omega^nB\TOP(n)) \ar[d] \ar[r] &  H^{n-1}(\Omega^n(\TOP/\TOP(n))) \\
H^{n}(\Omega^{n}B\TOP) \ar[d] \\
H^{n}(\Omega^{n}B\TOP(n))
}
\end{split}
\end{equation}
The horizontal arrow is induced by the canonical map from $\Sigma\Omega^{n}(\TOP/\TOP(n))$ to the
mapping cone of the inclusion $\Omega^{n}B\TOP(n)\to \Omega^{n}B\TOP$. The column is a part of the
long exact sequence of the pair $(\Omega^{n}B\TOP,\Omega^{n}B\TOP(n))$.
By hypothesis~\ref{hypa}, the class
\[ \Omega^{n}p_{n/2}\in H^{n}(\Omega^{n+2}B\TOP) \]
maps to zero in
$H^{n}(\Omega^{n}B\TOP(n))=H^{n}(\Omega^{n}B\STOP(n))$, since $\Omega^{n}p_{n/2}=\Omega^{n}p_{n/2}-\Omega^{n}e^2$.
Therefore $\Omega^{n}p_{n/2}$ lifts to a
class in $H^{n}(\Omega^{n}B\TOP,\Omega^{n}B\TOP(n))$. Let $z'$ be the image of that
in $H^{n-1}(\Omega^{n}(\TOP/\TOP(n)))$.
It is easy to show that $z'$ maps to a generator of
\[ H^{n-3}_{S^1}(\Omega^{n+2}Y_\QQ)\cong\QQ\,,\]
which we view as a direct summand of $H^{n-3}_{S^1}(\Omega^{n+2}Y)$,
by comparing diagram~(\ref{eqn-AtoB}) with a similar diagram where $\TOP$ and $\TOP(n)$ are replaced by $\Or$
and $\Or(n)$.
\qed

\bigskip
We finish this section with a very short preview of the rest of the paper. Since we have already obtained all of
\thehypa $\Rightarrow$ \thehypb and almost all of \thehypc $\Rightarrow$ \thehypa, the main point is to establish
\thehypb $\Rightarrow$ \thehypc. This is hard.


\section{Orthogonal calculus}\label{sec-ortho}
Let $\Ofun$, $\Tfun$ and $\Gfun$
be the continuous functors on $\sJ$ given for $V$ in $\sJ$ by $V\mapsto B\Or(V)$,
$V\mapsto B\TOP(V)$ and $V\mapsto B\G(V)$, respectively, where $\G(V)$ is the
topological group-like monoid of homotopy equivalences $S(V)\to S(V)$.
By orthogonal calculus \cite{WeissOrtCal}, the
functors $\Ofun$, $\Tfun$ and $\Gfun$ determine spectra
\[ \Theta \Ofun^{(i)},\quad\Theta
\Tfun^{(i)},\quad \Theta \Gfun^{(i)} \]
with an action of $\Or(i)$, for any integer $i>0$.
These are the $i$-th derivatives at infinity of $\Ofun$, $\Tfun$ and $\Gfun$, respectively.
The inclusions
\[ \Ofun \lra \Tfun \lra \Gfun
\]
determine maps of spectra
\begin{equation}\label{eqn-mapofspectra}
\Theta\Ofun^{(i)}\lra \Theta\Tfun^{(i)} \lra \Theta\Gfun^{(i)}
\end{equation}
which respect the actions of $\Or(i)$. \newline
\emph{Aside} (i): strictly speaking there
is no inclusion $\Tfun\to\Gfun$. Instead there is a diagram of continuous
functors and natural transformations of the form
\[
\xymatrix{ \Tfun & \Gfun \ar[d]^\simeq  \\
\Tfun' \ar[u]^\simeq \ar[r] & \Gfun'
}
\]
where $\Tfun'(V)$ is the classifying space of the topological group of homeomorphisms $V\to V$ fixing $0$, 
and $\Gfun'(V)$ is the grouplike topological monoid of homotopy equivalences from $V\smin\{0\}$ to
itself. \newline
\emph{Aside} (ii): It is known \cite{Cernavskii, EdwardsKirby} that $\TOP(V)$ is a locally contractible
topological group.
This is not directly important here, though, because we are mainly interested in the weak homotopy type of
$B\TOP(V)$. The singular simplicial set of $\TOP(V)$ is a simplicial group $\TOP_s(V)$.
A $k$-simplex in $\TOP_s(V)$ is a homeomorphism from $\Delta^k\times V$ to
$\Delta^k\times V$ over $\Delta^k$.
The canonical map $B|\TOP_s(V)|\to B\TOP(V)$ is a weak homotopy equivalence by general model category
principles. Therefore the weak homotopy type of $B\TOP(V)$ is at least geometrically intelligible.
Note that $V\mapsto B\TOP_s(V)$ is not a continuous functor on $\sJ$, which is why we prefer to work with 
$V\mapsto B\TOP(V)$.

\medskip
The Taylor tower of $\Ofun$ consists of
approximations $\Ofun\to T_i\Ofun$ for every $i\ge 0$, and maps $T_i\Ofun\to T_{i-1}\Ofun$ \
for $i>0$ under $\Ofun$. The homotopy fiber $L_i\Ofun$ of the map
$T_i\Ofun \to T_{i-1}\Ofun$ can be described as
\[  L_i\Ofun(V) ~\simeq~ \Omega^\infty
\Big(\big((V\otimes\RR^i)^c\wedge \Theta \Ofun^{(i)}\big)_{\ho \Or(i)}\Big) \]
(by a chain of natural homotopy equivalences). The functor $T_0\Ofun$ is essentially
constant, namely, $T_0\Ofun(V) \simeq B\Or$
by a chain of natural homotopy equivalences. The natural transformation $\Ofun\to T_i\Ofun$ has a
universal property, in the \emph{initial} sense: it is the best approximation of $\Ofun$ from the right by
a polynomial functor of degree $\le i$. \newline
Similarly, $\Tfun$ and $\Gfun$ have a
Taylor tower whose layers $L_i\Tfun$ and $L_i\Gfun$ for $i>0$ are determined, up to natural weak equivalence, by
the $i$-th derivative spectra of $\Tfun$ and $\Gfun$, respectively. Also, $T_0\Tfun$ is
essentially constant with value $B\TOP$ and $T_0\Gfun$ is essentially constant
with value $B\G$. There is one aspect in which $\Tfun$ differs substantially from
$\Ofun$ and $\Gfun$~: the Taylor towers of $\Ofun$ and $\Gfun$ are known \cite{Arone} to converge
to $\Ofun$ and $\Gfun$ respectively, that is,
\[ \Ofun(V) \simeq \holimsub{i} T_i\Ofun(V)~, \qquad \Gfun \simeq \holimsub{i} T_i\Gfun(V)~. \]
It is not known whether this holds for $\Tfun$, and Igusa's work on concordance stability
\cite{Igusa} indicates that it will not be easy to decide. \newline
This chapter analyses  the Taylor towers of the three functors $\Ofun$, $\Tfun$ and $\Gfun$ up to stage $2$ at most,
concentrating on rational aspects. This is done in part for illustration of methods. We are quite
aware that Arone~\cite{Arone} has already given an exhaustive, integral and very pretty description of the
Taylor tower of $\Ofun$. Our methods are more elementary and admittedly more pedestrian.

\medskip
We begin with the orthogonal calculus analysis of the functor $\Sfun\co \sJ \to {\bf T}_\pt$ given by
$V\mapsto V^c$, where $V^c$ is the one-point compactification
of the vector space $V$. Let $\Sigma_n$ be the symmetric group in $n$ letters.

\begin{prop}\label{prop-Ccalc}
The functor $\Sfun$ is rationally polynomial of degree 2, except for a possible deviation at $V=0$.
The first and second derivative spectra are
\[ \Theta \Sfun^{(j)}\simeq (\Or(j)/\Sigma_j)_+\wedge\Omega^{j-1}\underline{S}^0~,\]
for $j=1,2$, where the $\Or(j)$-action is trivial on the sphere spectrum and is the usual action on  $\Or(j)/\Sigma_j$.
\end{prop}

\proof Our proof is based on a relationship (probably first noted by Goodwillie) between orthogonal calculus and Goodwillie's calculus of
homotopy functors. Let $F$ be any homotopy functor from based CW-spaces to based spaces. Although
we do not routinely assume that $F$ is continuous, there is a standard construction which turns it into
a continuous functor. Assuming that $F$ is continuous, the functor $F\circ \Sfun$ on $\sJ$ is also continuous
and we can apply orthogonal calculus to it. \newline
Suppose for a start that $F$ is
\emph{homogeneous} of degree $n$. Then up to a chain of natural weak equivalences, $F$ has the form
\[  X \mapsto \Omega^\infty(X^{(n)}\wedge_{\Sigma_n}\mathbf \partial_nF) \]
where $\partial_nF$ is a CW-spectrum with a free (away from the base point) cellular action of $\Sigma_n$,
and $X^{(n)}$ is the $n$-fold smash power. It follows that in orthogonal calculus, the functor
$V\mapsto F(V^c)$ is also homogeneous of degree $n$. Namely, we can write
\[
\begin{array}{rcl}
F(V^c) & \simeq & \Omega^\infty((V\otimes\RR^n)^c\wedge_{\Sigma_n}\partial_nF) \\
& \cong & \Omega^\infty((V\otimes\RR^n)^c\wedge_{\Or(n)}(\Or(n)_+\wedge_{\Sigma_n}\partial_nF))~.
\end{array}
\]
Now let $F$ be the identity functor from based spaces to based spaces. From Goodwillie calculus,
we have polynomial approximations $F\to P_iF$ and the homogeneous layers
$\Lambda_iF$, so that there are natural homotopy fiber sequences
\[  \Lambda_iF\to P_iF\to P_{i-1}F \]
for $i>0$. We show by induction on $i$ that the functor $P_iF\circ \Sfun$ on $\sJ$ is almost polynomial of degree $\le i$,
in the sense that $(P_iF\circ \Sfun)(V)\to T_i(P_iF\circ \Sfun)(V)$ is an equivalence for all $V\ne 0$.
The induction beginning, $i=0$, is trivial. In fact the case $i=1$ is also clear because $P_1F=\Lambda_1F$ is
homogeneous of degree $1$.
For the induction step we assume $i>1$ and we have a homotopy fiber sequence
\[  \Lambda_iF\circ\Sfun\lra P_iF\circ \Sfun \to P_{i-1}F\circ\Sfun \]
where $\Lambda_iF\circ\Sfun$ is homogeneous of degree $i$ and $P_{i-1}F\circ\Sfun$ is polynomial of degree $\le i-1$.
The operator $T_i$ (the orthogonal calculus analogue of Goodwillie's $P_i$) respects homotopy fiber sequences,
so that we have a diagram
\[
\xymatrix{ \Lambda_iF\circ\Sfun \ar[r] \ar[d]^-\simeq & P_iF\circ \Sfun \ar[r] \ar[d] & P_{i-1}F\circ\Sfun \ar[d]^-{\simeq} \\
T_i(\Lambda_iF\circ\Sfun) \ar[r]  & T_i(P_iF\circ \Sfun) \ar[r]  & T_i(P_{i-1}F\circ\Sfun)
}
\]
We want to deduce that the middle column is also a homotopy equivalence, like the outer columns. This
is true when $V\ne 0$ because then the term $P_{i-1}F\circ \Sfun(V)=P_{i-1}F(V^c)$
is connected.
\newline
To continue we need to know something about the derivative spectra $\partial_nF$. It is known that $\partial_1F$ is a sphere
spectrum and that $\partial_2F$ is a looped sphere spectrum with a trivial action of $\Sigma_2$. This leads to the
formulae which we have given for the first and second derivative spectra of the functor $\Sfun=F\circ\Sfun$.
Finally we note that the natural map $F(X)\to P_2F(X)$ is a rational homotopy equivalence for $X=S^n$, assuming $n>0$.
(See remark~\ref{rem-Snaith} below.)
Therefore the natural map $\Sfun(V)=F(\Sfun(V))\to P_2F(\Sfun(V))$ is a rational homotopy equivalence when $V\ne 0$.
\qed

\begin{rem} \label{rem-Snaith} {\rm Let $X$ be a based space, $X^{\wedge k}$ the $k$-fold smash power, $Q(X)=\Omega^\infty\Sigma^\infty X$.
The symmetric group $\Sigma_k$ acts on $X^{\wedge k}$. The Snaith splitting of $\Omega^k\Sigma^kX$
means for $k=\infty$ that there is a zigzag of natural weak equivalences relating $\Sigma^\infty Q(X)$ to
\[ \bigvee_{k\ge 1} \Sigma^\infty ((X^{\wedge k})_{h\Sigma_k}). \]
It implies (modulo that zigzag) a natural map from $\Sigma^\infty(Q(X)/X)$ to $\Sigma^\infty((X^{\wedge 2})_{h\ZZ/2})$
which is $(3c+2)$-connected if $X$ is $c$-connected, $c>0$. The adjoint of that is a natural map from
$Q(X)/X$ to  $Q((X^{\wedge 2})_{h\ZZ/2})$
which is still $(3c+2)$-connected. Writing $j$ for the
composition of $Q(X)\to Q(X)/X$ with that last map, we get an inclusion
\[  X\lra \hofiber[j\co Q(X)\to Q((X^{\wedge 2})_{h\ZZ/2})] \]
which is still $3c$-connected. Since the functor $X\mapsto \hofiber[j\co Q(X)\to Q((X^{\wedge 2})_{h\ZZ/2})]$
is clearly polynomial of degree $2$, and since it comes with a natural transformation from $F=\id$ with these good approximation
properties, that natural transformation is the second Taylor approximation $F\to P_2F$. With this description
of $P_2F$, it is clear that $P_2F$ of a sphere has the same rational homotopy groups as the sphere, and since
the approximation $X\to P_2F(X)$ is so highly connected in the case of a sphere, it must be a rational equivalence.
}
\end{rem}

\medskip
We now generalize the functor $\Sfun$ of the previous proposition to allow wedge sums of shifts of $\Sfun$.
Let $J=(k_j)_{j=1,2,\dots}$ be a (finite or infinite) sequence of nonnegative integers. Then we
define $\Sfun_J$ to be the functor
\begin{equation}\label{def-generalR}
 V\mapsto \bigvee_j \Sfun(V\oplus\RR^{k_j})~.
\end{equation}

 \begin{prop}\label{prop-Rcalc}
The first and second derivative spectra of $\Sfun_J$ are as follows:
\[ \Theta \Sfun_J^{(1)}= \bigvee_j \Or(1)_+\wedge \underline{S}^{k_j}~, \]
\[ \Theta \Sfun_J^{(2)} = \Big(\bigvee_j  \Or(2)_+\wedge_{\Sigma_2}\Omega\underline{S}^{2k_j}\Big)
~\vee~ \Big(\bigvee_{j<\ell} \Or(2)_+\wedge \Omega\underline{S}^{k_j+k_\ell}\Big)~, \]
where $\Or(i)$ acts by translation on the first smash factor (and the symmetric group $\Sigma_2$ acts on $\Omega\underline{S}^{2k_j}=
\Omega(S^{k_j}\wedge S^{k_j}\wedge\underline{S}^0)$ by permuting the two copies of $S^{k_j}$).
\end{prop}

\proof As in the previous proof let $F$ be the identity functor on based spaces.
We try $P_iF\circ \Sfun_J$ as a candidate for $T_i\Sfun_J$. First we need to show that
$\Lambda_iF\circ \Sfun_J$ is homogeneous of degree $i$ for $i>0$. There is the standard
expression
\[
\begin{array}{rcl}
(\Lambda_iF\circ \Sfun_J)(V) & = & \Lambda_iF\big(\bigvee_j V^c\wedge S^{k_j}\big) \\
& \simeq &\!\!
\Omega^\infty\big(\big(\bigvee_{(j_1,\dots,j_i)}(V\otimes\RR^i)^c\wedge S^{\sum k_{j_i}}\big)
\wedge_{\Or(i)}\big(\Or(i)_+\wedge_{\Sigma_i}\partial_iF\big)\big) \\
& = &\!\!\Omega^\infty\big(\bigvee_{(j_1,\dots,j_i)} (V\otimes\RR^i)^c
\wedge_{\Or(i)}\big(S^{\sum k_j}\wedge \Or(i)_+\wedge_{\Sigma_i}\partial_iF\big)\big) \\
& = &\!\! \Omega^\infty\big((V\otimes\RR^i)^c
\wedge_{\Or(i)}\big(\bigvee_{(j_1,\dots,j_i)} S^{\sum k_j}
\wedge \Or(i)_+\wedge_{\Sigma_i}\partial_iF\big)\big)
\end{array}
\]
which shows that $\Lambda_iF\circ \Sfun_J$ is homogeneous of degree $i$. We deduce as before that
$P_iF\circ\Sfun_J$ is almost
polynomial of degree $i$, that is, the natural map
\[ (P_iF\circ\Sfun_J)(V)\to T_i(P_iF\circ\Sfun_J)(V) \]
is a homotopy equivalence for $V\ne 0$.
Also the natural map
\[ \Sfun_J(V)= F(\Sfun_J(V)) \lra P_iF(\Sfun_J(V))  \]
is $((i+1)(\dim(V)-1)-k)$-connected for a fixed $k$ independent of $V$ and $i$,
by the convergence of the Goodwillie tower for $F$. By the construction of $T_i$~, this implies that
\[ T_i\Sfun_J\lra T_i(P_iF\circ \Sfun_J) \]
is a weak equivalence. (Let $\Efun$ be a continuous functor from $\sJ$ to based spaces. Then $T_i\Efun$
is defined by iterating the construction $\tau_i$ defined by
\[  \tau_i\Efun(V) = \holimsub{0\ne U\le\RR^{i+1}} \Efun(U\oplus V). \]
See \cite{WeissOrtCal}. Here we are looking at a topological homotopy inverse limit, and $U$ runs over the nonzero linear
subspaces of $\RR^{i+1}$. If a natural transformation $\Efun_0\to \Efun_1$ has the
property that $\Efun_0(V)\to \Efun_1(V)$ is $((i+1)(\dim(V)-1)-k)$-connected for a fixed $k$ independent of $V$ and $i$,
then the induced natural map $\tau_i \Efun_0(V)\to \tau_i \Efun_1(V)$ is $((i+1)(\dim(V)-1)-k+1)$-connected.)
\newline
Now we have a chain of natural weak equivalences relating $T_i\Sfun_J(V)$ to $(P_iF\circ\Sfun_J)(V)$ for $V\ne 0$.
Therefore the homogeneous layers in the Taylor tower of $\Sfun_J$ are the $\Lambda_iF\circ\Sfun_J$. The description
of the derivative spectra of $\Sfun_J$ for $i=1,2$ follows as in the previous proof. \qed

\medskip
Now we shall sketch the orthogonal calculus analysis of $\Ofun$, concentrating
on rational aspects where that saves energy. Much more detailed results can be found in~\cite{Arone}.

\begin{prop}\label{prop-calcE}
The functor $\Ofun$ is rationally polynomial of degree 2, in the sense that the canonical map
$\Ofun(V)\to T_2\Ofun(V)$ is a rational homotopy equivalence over $B\ZZ/2$,
for every $V\ne 0$. The derivative spectra of $\Ofun$ are as follows
\begin{itemize}
 \item[(i)] $\Theta \Ofun^{(1)}\simeq \underline{S}^0$ with trivial action of $\Or(1)$;
 \item[(ii)] $\Theta \Ofun^{(2)}\simeq \Omega\underline{S}^0$ with rationally trivial action of $\Or(2)$.
\end{itemize}
\end{prop}

\proof
By definition the spectrum $\Theta \Ofun^{(1)}$
is made up of the based spaces
\[ \Theta \Ofun^{(1)}(n)=\hofiber[ \Ofun(\RR^n) \to \Ofun(\RR^{n+1})]~~\simeq~S^n \]
and so turns out to be a sphere spectrum $\underline{S}^0$. The generator of $\Or(1)$
acts on $\Theta \Ofun^{(1)}(n)$ alias $S^n=\RR^n\cup\{\infty\}$ via $-\id\co \RR^n\to \RR^n$. The structure maps
\[  S^1\wedge \Theta \Ofun^{(1)}(n) \to \Theta \Ofun^{(1)}(n+1)  \]
are $\Or(1)$-maps, where we use the standard conjugation action on $S^1$ and the
resulting diagonal action on $S^1\wedge \Theta \Ofun^{(1)}(n)$. Therefore, strictly speaking, the
structure maps are in a twisted relationship to the actions of $\Or(1)$ on the various $\Theta \Ofun^{(1)}(n)$,
but there are mechanical ways to untwist this (see also the beginning of section~\ref{sec-orcalsmooth} below)
and the result is a sphere spectrum with trivial action of $\Or(1)$. \newline
For the description of the second derivative spectrum we reduce this to proposition~\ref{prop-Ccalc}.
We have a natural homotopy fiber sequence
\[  \Sfun(V)\lra \Ofun(V)\lra \Ofun(V\oplus \RR)  \]
inducing a corresponding homotopy fiber sequence of spectra
\[ \Theta \Sfun^{(2)}\lra  \Theta \Ofun^{(2)}\wedge S^0 \lra  \Theta \Ofun^{(2)}\wedge S^2~, \]
where we think of $S^2$ as $(\RR^2)^c$ and the maps of this homotopy fiber sequence preserve the $\Or(2)$ actions
(in particular, the $\Or(2)$-action on $\Theta \Ofun^{(2)}\wedge S^2$ is the diagonal one).
Consequently, we have that
\begin{equation}\label{eqn-theta2c}
   \Theta \Sfun^{(2)} \simeq \Omega( \Theta \Ofun^{(2)}\wedge (S^2/S^0))\simeq \Theta \Ofun^{(2)}\wedge S^1_+.
\end{equation}
Taking homotopy orbit spectra for the action of $\SOr(2)$ and using our previous formula for $\Theta \Sfun^{(2)}$
we obtain that
\[  \Omega\underline{S}^0 \simeq \Theta \Ofun^{(2)}~.\]
This equivalence does not fully keep track of
$\Or(2)$ actions, but it does allow us to say that orientation reversing elements of $\Or(2)$ act by self-maps homotopic
to the identity.  Consequently the action of $\Or(2)$ on $\Theta \Ofun^{(2)}$ is rationally trivial, in the sense
that $\Theta \Ofun^{(2)}$ is rationally weakly equivalent as an $\Or(2)$-spectrum to an $\Or(2)$-spectrum with
trivial action of $\Or(2)$. To see that $\Ofun$ is rationally polynomial of degree 2 we consider the commutative diagram
\[  \xymatrix{ \Sfun(V) \ar[r]\ar[d] & \Ofun(V) \ar[r]\ar[d] & \Ofun(V\oplus \RR) \ar[d] \\
T_2\Sfun(V) \ar[r] & T_2\Ofun(V) \ar[r] & T_2\Ofun(V\oplus \RR)~,  } \]
where the rows are homotopy fiber sequences. If $V\ne 0$ the left-hand vertical arrow is a rational equivalence
and so the right-hand square is rationally
a homotopy pullback square over $B\ZZ/2$. Therefore, by iteration,
\[  \xymatrix{  \Ofun(V) \ar[r]\ar[d] & \Ofun(V\oplus \RR^\infty) \ar[d] \\
 T_2\Ofun(V) \ar[r] & T_2\Ofun(V\oplus \RR^\infty)~,  } \]
is rationally a homotopy pullback square over $B\ZZ/2$. Here the right-hand column is a
homotopy equivalence. So the left-hand column is a rational homotopy equivalence over $B\ZZ/2$.
\qed

\bigskip
The space $S\G(\RR^n)=S\G(n)$ can be investigated using the homotopy fiber sequence
\begin{equation} \label{eqn-SG}
\CD \Omega_1^{n-1}S^{n-1} @>>> S\G(n) @>>> S^{n-1} \endCD
\end{equation}
where the right-hand map is evaluation at the base point of $S^{n-1}$. For even
$n>0$ we can deduce immediately
\[ S\G(n)\simeq_{\QQ} S^{n-1}\quad,\quad BS\G(n) \simeq_{\QQ} K(\QQ,n).~ \]
For odd $n>1$ the connecting homomorphisms $\pi_{n-1}S^{n-1}\to \pi_{n-2}\Omega_1^{n-1}S^{n-1}$ in the
long exact homotopy group sequence of~(\ref{eqn-SG}) are rational isomorphisms, as can be seen
by comparing~(\ref{eqn-SG}) with the homotopy fiber sequence
\[ \SOr(n-1)\to \SOr(n) \to S^{n-1}~. \]
Therefore
\[ S\G(n)\simeq_{\QQ} K(\QQ,2n-3)\quad,\quad BS\G(n) \simeq_{\QQ} K(\QQ,2n-2)~. \]
For $B\G(n)$ with arbitrary $n\ge 2$ we obtain therefore
\begin{equation} \label{eqn-BGn}
B\G(n) \simeq_{\QQ} \left\{ \begin{array}{ll}
K(\QQ,n)_{\ho \ZZ/2} & \textup{ $n$ even, $n>0$} \\
K(\QQ,2n-2) & \textup{ $n$ odd, $n>1$}
\end{array}
\right.
\end{equation}
where $\ZZ/2$ acts by sign change on the $\QQ$ in $K(\QQ,n)$. (In the case $n$ even, this is a rational
homotopy equivalence over $B\ZZ/2$.)
Furthermore, it follows from~(\ref{eqn-SG}) and~(\ref{eqn-BGn}) that
for odd $n>1$, the diagram
\begin{equation} \label{eqn-BGhofi} 
\xymatrix@C=40pt{
S\G(n-1) \ar[r]^{\textup{inc.}} &  S\G(n) \ar[r]^{\textup{eval. at $\pt$}} & S^{n-1}
}
\end{equation}
is a rational homotopy fiber sequence.
These calculations can be summarized as follows. In the
case of even $n$, the twisted Euler class in
$H^n(B\G(n);\ZZ^t)$ (with local coefficients $\ZZ^t$ determined
by the nontrivial action of $\pi_1B\G(n)\cong\ZZ/2$ on $\ZZ$)
detects the entire rational homotopy of $B\G(n)$. In the case of
odd $n>1$, there is a class in $H^{2n-2}(B\G(n);\QQ)$ which detects the entire rational homotopy of $B\G(n)$,
and extends the squared Euler class $e^2\in H^{2n-2}(BS\G(n-1);\QQ)$.

\begin{prop}\label{prop-calcxi}
The functor $\Gfun$ is rationally polynomial of degree 2, in the sense that the canonical map
$\Gfun(V)\to T_2\Gfun(V)$ is a rational homotopy equivalence over $B\ZZ/2$,
for every $V$ of dimension $\ge 2$.
The derivative spectra of $\Gfun$ satisfy the following:
\begin{itemize}
 \item[(i)] the map $\Theta \Ofun^{(1)}\to\Theta \Gfun^{(1)}$ induced by $\Ofun\hookrightarrow\Gfun$
 is a homotopy equivalence;
 \item[(ii)] the map $\Theta \Ofun^{(2)}\to\Theta \Gfun^{(2)}$ induced by
$\Ofun\hookrightarrow\Gfun$ fits into a commutative triangle of spectra
with $\Or(2)$ action
\[
\xymatrix@R=17pt@C=30pt{\Theta \Ofun^{(2)} \ar[d] \ar[rd] &  \\
\Theta \Gfun^{(2)} \ar[r]_-{\simeq_\QQ} &  \map(S^1,\Theta\Ofun^{(2)})
}
\]
where $\Or(2)$ acts in the standard manner on $S^1$ and the vertical arrow is given by inclusion of the constant maps.
\end{itemize}
\end{prop}

\proof We begin by showing (i) and (ii).
It is known~\cite{Haef,Wall} that the map from $\Or(n+1)/\Or(n)$ to $\G(n+1)/\G(n)$, induced by inclusion, is $(2n-c)$-connected
(for a small constant $c$
independent of $n$). It follows
immediately that the natural map of first derivative spectra induced by $\Ofun\to\Gfun$
is a weak homotopy equivalence. Therefore the first derivative spectrum of $\Gfun$ is (equivalent to) a
sphere spectrum with trivial action of $\Or(1)$.
For the second derivative spectrum we use the homotopy fiber sequence
\begin{equation}\label{eqn-Xicalc2}
 \Omega_{\pm 1}^V\Sfun(V) \lra \Omega \Gfun(V\oplus \RR)\lra \Sfun(V)~,
\end{equation}
where $\Omega^V(X)$ is the space of pointed maps from $V^c$ to $X$ (for $X$ a based space) and the $\pm 1$ singles out
the degree $\pm 1$ components.
In lemma~\ref{lem-OmegaXi} below we show that the functor taking $V$ to $\Omega_{\pm 1}^V\Sfun(V) $ is rationally
polynomial of degree 1, for all nonzero $V$.
Therefore the homotopy fiber sequence (\ref{eqn-Xicalc2}) induces a rational equivalence of second derivative spectra,
\[  \Omega(\Theta \Gfun^{(2)}\wedge S^2)\simeq_\QQ \Theta \Sfun^{(2)}~. \]
(See lemma~\ref{lem-specshift}.) More precisely, we have a commutative diagram of second derivative spectra
\begin{equation} \label{eqn-comparebgboc}
\begin{split}
\xymatrix@R=15pt@C=35pt{\Omega(\Theta \Ofun^{(2)}\wedge S^2) \ar[d]\ar[rd]^{\textup{eval.}} &\\
  \Omega(\Theta \Gfun^{(2)}\wedge S^2)\ar[r]^{\quad\textup{eval.}} & \Theta \Sfun^{(2)}~,  }
\end{split}
\end{equation}
where the horizontal map is a rational homotopy equivalence. We saw in the proof of proposition~\ref{prop-calcE}  that the diagonal arrow
in the diagram is in fact the map
\begin{equation}\label{eqn-loopXi}
\Omega(\Theta \Ofun^{(2)}\wedge S^2)
~~\lra~~\Omega(\Theta \Ofun^{(2)}\wedge (S^2/S^0))\simeq\Theta \Ofun^{(2)}\wedge S^1_+ ~,
\end{equation}
where we use (\ref{eqn-theta2c}).
Therefore we can change triangle (\ref{eqn-comparebgboc}) into a commutative triangle
\begin{equation} \label{eqn-comparebgboc2}
\begin{split}
\xymatrix@R=15pt@C=35pt{\Omega(\Theta \Ofun^{(2)}\wedge S^2)
\ar[d]\ar[rd]^{(\ref{eqn-loopXi})} &\\
\Omega(\Theta \Gfun^{(2)}\wedge S^2)\ar[r]_{\simeq_\QQ} & \Theta \Ofun^{(2)}\wedge S^1_+ ~.  }
\end{split}
\end{equation}
By undoing the looping and the double suspension, we obtain the commutative triangle
\begin{equation} \label{eqn-comparebgboc3}
\begin{split}
\xymatrix@R=15pt@C=35pt{\Theta \Ofun^{(2)} \ar[d]\ar[rd] &\\
  \Theta \Gfun^{(2)}\ar[r]_-{\simeq_\QQ} & \Omega^{2}(\Theta \Ofun^{(2)}\wedge (S^2/S^0)) ~.  }
\end{split}
\end{equation}
There is a homotopy equivalence from $\Omega^{2}(\Theta \Ofun^{(2)}\wedge (S^2/S^0))$
to  $\map(S^1,\Theta \Ofun^{(2)})$, preserving $\Or(2)$-actions,  which we describe in adjoint form by
\[ S^1_+\wedge \Omega^{2}(\Theta \Ofun^{(2)}\wedge (S^2/S^0)) \lra  \Theta \Ofun^{(2)}.\]
Namely, every choice of point $z$ in $S^1$ determines a nullhomotopy for the inclusion of $S^0$ in $S^2$ and thereby
a map $S^2/S^0 \to S^2\vee S^1 \to S^2$. So we have
\[ z_+\wedge \Omega^{2}(\Theta \Ofun^{(2)}\wedge (S^2/S^0)) \lra \Omega^{2}(\Theta \Ofun^{(2)}\wedge S^2)\simeq \Theta \Ofun^{(2)}\]
for every $z\in S^1$, and using these for all $z$ gives the required map. It is equivariant for the diagonal
$\Or(2)$-action on the source. \newline 
The homotopy fiber sequence~(\ref{eqn-Xicalc2}) implies that the approximation
\[ \Gfun(V)\to T_2\Gfun(V) \]
is a rational
homotopy equivalence over $B\ZZ/2$ when $\dim(V)\ge 2$.

\begin{lem}\label{lem-OmegaXi}
The functor $V\mapsto \Omega_{\pm 1}^V\Sfun(V)$ is rationally polynomial of degree 1, for nonzero $V$.
\end{lem}

\proof
We want to show that
for nonzero $V$ the approximation
\[ \Omega_{\pm 1}^V\Sfun(V)\to T_1(\Omega_{\pm 1}^V\Sfun)(V) \]
is a rational
homotopy equivalence of componentwise nilpotent spaces.
Let $\Efun_u(V)=\Omega_{\pm 1}^V\Sfun(V)$ and let $\Efun$ be the slightly simpler functor defined by
$\Efun(V)=\Omega^V\Sfun(V)$.
It is enough to show that $\Efun$ is rationally polynomial of degree $\le 1$, allowing a deviation at $V=0$,
because there is a homotopy pullback square
\[
\CD
\Efun_u(V) @>>> \Efun(V) \\
@VVV  @VVV \\
\{\pm 1\} @>>> \ZZ
\endCD
\]
(for $V\ne 0$) where the functors in the lower row are of degree $0$. \newline
The functors $\Efun$ and $\Ofun$ are closely related by the following construction. Given any continuous
functor $\Ffun$ from $\sJ$ to based spaces, the continuous functor $\Ffun^{(1)}$ defined by
$\Ffun^{(1)}(V):= \hofiber[\Ffun(V)\to \Ffun(V\oplus\RR)]$
has additional structure in the shape of binatural maps
\[  V^c\wedge \Ffun^{(1)}(W) \lra \Ffun^{(1)}(V\oplus W). \]
These maps enjoy associativity and unital properties, making $\Ffun^{(1)}$ into what is nowadays called
an orthogonal spectrum. Consequently $V\mapsto \Omega^V\Ffun^{(1)}(V)$ is again a continuous functor
from $\sJ$ to spaces; and it is the functor $\Efun$
if we start out with $\Ffun=\Ofun$.
Now we continue with some general observations.
\begin{itemize}
\item[(i)] \emph{If $\Ffun$ is homogeneous of degree $n$, then $V\mapsto \Omega^V\Ffun^{(1)}(V)$ is
homogeneous of degree $n-1$.}
This is shown in \cite[Expl.~5.7]{WeissOrtCal}.
\item[(ii)] \emph{If $\Ffun$ is polynomial of degree $\le n$, then
$V\mapsto \Omega^V\Ffun^{(1)}(V)$ is polynomial of degree $\le n-1$.}
This follows easily from (i) by writing $\Ffun$ as a finite tower with homogeneous layers.
\end{itemize}
Since we want to
apply a rational version of (ii), relying on proposition~\ref{prop-calcE}, we need to be aware of $\pi_0$ related problems.
Then again, since $\pi_0\Efun(V)$ for $V\ne 0$ is what it is, there are no such problems. \qed

\begin{rem} \label{rem-goodextzed} Let $\Gfun^u(V)=\hofiber[\Gfun(V)\to T_0\Gfun(V)]$.
Then clearly $L_i\Gfun^u=L_i\Gfun$ for
$i>0$. From the rational homotopy fiber sequence (\ref{eqn-BGhofi}) we obtain, for odd $n>2$, another rational
homotopy fiber sequence
\[ S^{n-1}\lra \Gfun^u(\RR^{n-1}) \lra  \Gfun^u(\RR^{n})~.\]
It follows that the inclusion-induced map $\Gfun^u(\RR^{n-1}) \lra  \Gfun^u(\RR^{n})$ is not rationally
nullhomotopic. (Indeed, $\Gfun^u(\RR^{n-1})$ is not a rational homotopy retract of $S^{n-1}$, since
the square of the Euler class in the rational cohomology of $\Gfun^u(\RR^{n-1})$ is nonzero.) \newline
By proposition~\ref{prop-calcxi}, the above rational homotopy fiber sequence can also be written in the form
$S^{n-1}\to T_2\Gfun^u(\RR^{n-1}) \to  T_2\Gfun^u(\RR^{n})$.
As $L_2\Gfun(V)$ is rationally trivial for even-dimensional $V\ne 0$ (see also remark~\ref{rem-pd}
below) and $L_1\Gfun$ is
rationally trivial for odd-dimensional $V$, the nonexistence
of a nullhomotopy for $T_2\Gfun^u(\RR^{n-1})\to T_2\Gfun^u(\RR^{n})$, for odd $n>2$, implies
that the natural homotopy fiber sequence
\begin{equation} \label{eqn-amazingnonsplit} L_2\Gfun(V) \lra T_2\Gfun^u(V) \lra L_1\Gfun(V) \end{equation}
does not admit a natural rational splitting. To say this more carefully, it is not
possible to produce another continuous functor $\Efun$ from $\sJ$ to spaces and a
natural transformation $T_2\Gfun^u\to \Efun$ such that the resulting composition
\[ L_2\Gfun(V)\lra T_2\Gfun^u(V)\lra \Efun(V) \]
is a rational equivalence for every $V$ of dimension $\gg0$.
Note that~(\ref{eqn-amazingnonsplit}) does admit a rational splitting for every $V\ne 0$ individually,
because one of $L_2\Gfun(V)$, $L_1\Gfun(V)$ is rationally contractible.
\end{rem}

\begin{rem}\label{rem-pd} Let $M$ be a closed smooth manifold
and $\underline{E}$ a CW-spectrum.
A Poincare duality principle identifies $\map(M,\underline{E})$ with
\[ \underline{F}= \int_{x\in M} \Omega^{T_xM}\underline{E}. \]
Here the spectra $\Omega^{T_xM}\underline{E}$ for $x\in M$ together make up
a fibered spectrum over $M$, and $\underline{F}$ is the ordinary spectrum
obtained by passing to total spaces and collapsing the zero sections. Such an
identification can also be used when $M$ and $\underline{E}$ come with actions
of a compact Lie group $G$. In particular, for $G=\Or(2)$ and $M=S^1$ with the standard
action of $\Or(2)$ and $\underline{E}=\Omega\underline{S}^0$
with the trivial action, the spectrum $\map(M,\underline{E})$ can be described as
$\Omega^2(S^1_+\wedge\underline{S}^0)$
with the following action of $\Or(2)$: trivial on the $\underline{S}^0$ factor, standard on the $S^1_+$
factor, \emph{adjoint} action on one of the loop coordinates
(the action of $\Or(2)$ on its Lie algebra).
\end{rem}

\bigskip

We turn to the (rational) Taylor tower of the functor $\Tfun$.
The inclusion $\Ofun \to \Tfun$ induces a rational homotopy
equivalence $T_0\Ofun\to T_0\Tfun$, which just restates the Thom-Novikov result
\[ B\Or\simeq_\QQ B\TOP~.\]
The spectrum $\Theta \Tfun^{(1)}$
has a known rational description with action of $\Or(1)$, due to Waldhausen, Borel and Farrell-Hsiang
\cite{Wald,Borel,FarrellHsiang}. Below, $\underline{A}(\pt)$ is Waldhausen's $A$-theory spectrum, also known as the
algebraic K-theory spectrum of the ring spectrum
$\underline{S}^0$, and $\underline{K}(\QQ)$ is the algebraic K-theory of the ring $\QQ$.

\begin{prop}\label{prop-calcF}
The functor $\Tfun$ has first derivative spectrum
\[ \Theta \Tfun^{(1)}\simeq \underline{A}(\pt) \simeq_\QQ \underline{K}(\QQ).\]
The $\Or(1)$ action is the standard duality action on $K$-theory. Hence
\[ \Theta \Tfun^{(1)}~\simeq_{\QQ}~ \Theta \Ofun^{(1)} \vee \bigvee_{i=1}^{\infty} S^{4i+1}\wedge\underline{H}\ZZ~ \]
where $\Or(1)$ acts on the summand $S^{4i+1}\wedge\underline{H}\ZZ$ through its standard action on $\ZZ$.
\end{prop}

\proof
By definition the spectrum $\Theta \Tfun^{(1)}$ is made up of spaces
\[ \Theta \Tfun^{(1)}(n)=\TOP(n+1)/\TOP(n) \] with structure maps
analogous to those of $\Theta \Ofun^{(1)}$. The identification of the spectrum 
\[ \{\TOP(n+1)/\TOP(n)\,|\,n\in \NN\} \] 
with $\underline{A}(\pt)$
comes from~\cite{Wald}. It relies on the smoothing theory description of spaces of smooth $h$-cobordisms
over $D^n$, as in example~\ref{exa-smoothhcobdisc}. Modulo that it is a central part of Waldhausen's development
of the algebraic $K$-theory tradition in $h$-cobordism theory, a tradition which started with the
$h$-cobordism and $s$-cobordism theorems \cite{Smale, Milnorhcob, Cohen, Kervaire}.
The identification of the canonical $\Or(1)$-action on the spectrum $\{\TOP(n+1)/\TOP(n)\,|\,n\in \NN\}$ with the $\ZZ/2$-action
on $\underline{A}(\pt)$  by (Spanier-Whitehead) duality is due to~\cite{Vogell}, again going through $h$-cobordism theory.
The inclusion-induced
map $\Theta^{(1)}\Ofun\to \Theta^{(1)}\Tfun$, alias $\underline{S}^0\to \underline{A}(\pt)$, admits an \emph{integral}
left inverse (with $\Or(1)$-equivariance, in the homotopy category of $\Or(1)$-spectra). The rational equivalence
$\underline{A}(\pt) \simeq_\QQ \underline{K}(\QQ)$ is a consequence of the rational equivalence between the sphere spectrum
(as a ring spectrum) and the Eilenberg-MacLane spectrum $H\QQ$.
The calculation of the rational homotopy groups of $\underline{K}(\QQ)$ follows from the calculation of the rational
cohomology groups of  $B\GL(\QQ)$, due to Borel~\cite{Borel}. The result is
\[  \pi_n(\underline{K}(\QQ))\otimes\QQ~\cong~\left\{
 \begin{array}{ll} \QQ & n=0 \\
\QQ & n=5,9,13,17,\dots \\
0 & \textup{otherwise.}
  \end{array}
\right.
\]
The action of $\Or(1)$ on $\pi_n(\underline{K}(\QQ))$ is trivial for $n=0$ and nontrivial (sign change)
for $n=5,9,13.\dots$~.  \qed

\begin{rem} Hypothesis~\thehypc implies
hypothesis~\thehypa, for all even $n$. {\rm Most of this has been shown already in proposition~\ref{prop-CtoA}
but there we had to sacrifice the case $n=4$. The arguments here are much the same, though.
We start by reformulating hypothesis~\ref{hypa} as a statement about $B\TOP(n)$ instead of
$B\STOP(n)$, where $n$ is even. The Euler class $e$ in $H^{n}(B\STOP(n);\QQ)$ comes from an Euler class $e_t$ in
$H^{n}(B\TOP(n);\QQ^t)$ where $\QQ^t$ is a ``twisted'' local coefficient system, the twist being
determined by the first Stiefel-Whitney class of the universal euclidean bundle on $B\TOP(n)$.
It is therefore more than enough to show that
\[ e_t^2=p_{n/2}\in H^{2n}(B\TOP(n);\QQ) \]
if hypothesis \thehypc holds. \newline
Assuming hypothesis \thehypc, we have a functor splitting
$\Ofun\to \Tfun \to \Dfun$ such that the composition
$\Ofun\to \Dfun$ is a rational homotopy equivalence (over $B\ZZ/2$). Therefore we can speak of the twisted Euler class
$e_t$ in $H^{n}(\Dfun(\RR^n);\QQ^t)$ and the Pontryagin class $p_{n/2}$ in $H^{2n}(\Dfun(\RR^n);\QQ)$. For these we have
\[ e_t^2=p_{n/2}~. \]
It is therefore enough to show that under the map $\Tfun\to\Dfun$, the class $p_{n/2}\in H^{2n}(\Dfun(\RR^{n});\QQ)$
pulls back to the Pontryagin class in $H^{2n}(\Tfun(\RR^n);\QQ)$ and the class $e_t\in H^{n}(\Dfun(\RR^{n});\QQ^t)$ pulls back
to the twisted Euler class in  $H^n(\Tfun(\RR^n);\QQ^t)$. For the Pontryagin classes this follows from the commutativity of the diagram
\[
\xymatrix@R=16pt{
\Tfun(\RR^n) \ar[r]\ar[d] & T_0\Tfun(\RR^n)\ar[d]^{\simeq_\QQ} \\
\Dfun(\RR^n) \ar[r] & T_0\Dfun(\RR^n)~.
}
\]
(The Pontryagin classes come from the right-hand column.) For the Euler classes it follows from the commutativity of
the diagram
\[
\xymatrix@R=16pt{
\Tfun(\RR^n) \ar[r]\ar[d] & T_1\Tfun(\RR^n)\ar[d]^{\simeq_\QQ} \\
\Dfun(\RR^n) \ar[r] & T_1\Dfun(\RR^n)~.
}
\]
(The Euler classes come from the right-hand column.) The right-hand vertical arrow is a rational homotopy
equivalence over $B\ZZ/2$ by proposition~\ref{prop-calcF}, since $n$ is even.} \qed
\end{rem}

\begin{rem} \label{rem-eutopo} Assuming hypothesis~{\rm\ref{hypa}} we have
$p_{n/2}=e^2\in H^{2n}(B\STOP(n);\QQ)$. Hence $p_{n/2}$ is zero in
$H^{2n}(B\STOP(m);\QQ)$ when $m<n$. The restriction homomorphism
\[ H^{2n}(B\TOP(m);\QQ)\lra H^{2n}(B\STOP(m);\QQ) \]
is injective, as $B\STOP(m)$ is homotopy equivalent to a double cover of
$B\TOP(m)$. Therefore hypothesis~{\rm\ref{hypa}} implies that $p_{n/2}=0$ in $H^{2n}(B\TOP(m);\QQ)$
whenever $m<n$. \qed

\end{rem}

\section{Natural transformations}\label{sec-homoext}
Denote by $\sE_0$ the category of continuous functors from $\sJ$ to ${\bf T}_\pt$
and their natural transformations.
We want to promote $\sE_0$ to a model category in the sense of Quillen \cite{Hirschhorn, Hovey}.
For that we need to specify subcategories of weak equivalences, fibrations and cofibrations.
The weak equivalences and fibrations are defined levelwise. That is, a morphism
$\Efun\to \Ffun$ in $\sE_0$ is a \emph{weak equivalence} if $\Efun(V)\to \Ffun(V)$
is a weak homotopy equivalence for every $V$ in $\sJ$, and a \emph{fibration} if
$\Efun(V)\to \Ffun(V)$ is a Serre fibration for every $V$ in $\sJ$.
This choice of weak equivalences and fibrations determines the cofibrations, as in any model
category structure. It is known to be a consistent choice, i.e.,
a model category structure on $\sE_0$ with fibrations and weak
equivalences as specified exists. See for example \cite[Lemma 6.1]{BarnesOman} and \cite[Thm. 6.5]{MMSS}.

\medskip
The CW-functors defined in \cite{WeissOrtCal}
are convenient examples of cofibrant objects in  $\sE_0$\,.
For every continuous $\Efun$ from $\sJ$ to based $k$-spaces there exist
a based CW-functor $\Xfun$ and a weak equivalence $\Xfun\to \Efun$.
(Very briefly, a CW-functor $\Xfun$ is a functor from $\sJ$ to spaces equipped
with a sequence of subfunctors $\Xfun^i$ for $i=0,1,2,,\dots$ such that $\Xfun^0$
is a coproduct of co-representable functors
$\mor(V_\alpha,-)$ and $\Xfun^i$ for $i>0$ is the pushout of a diagram
\[   \Xfun^{i-1} \longleftarrow \coprod_\beta \mor(V_\beta,-)\times S^{i-1}~
\hookrightarrow~\coprod_\beta \mor(V_\beta,-)\times D^i ~, \]
and finally, $\Xfun=\colim_i~\Xfun^i$. It is an exercise to show that $\Xfun(0)$ is a CW-space with
$i$-skeleton $\Xfun^i(0)$.
We say that $\Xfun$ is a \emph{based} CW-functor if in addition a $0$-cell in $\Xfun(0)$
has been selected. Note that this determines a base point in $\Xfun(V)$ for every $V$, since $0$ is the
initial object in $\sJ$. Showing that based CW-functors are cofibrant boils down to showing mainly that functors of the form
$W\mapsto \mor(V,W)_+$
are cofibrant; this is easy by the Yoneda lemma.) A relative version of this idea leads to a good supply
of easy-to-understand cofibrations in the category of continuous functors from $\sJ$ to based $k$-spaces. \newline
Next we note that $\sE_0$ is a \emph{simplicial} model category \cite[4.2.18]{Hovey}
by means of an obvious action of the category of simplicial sets, namely
\[  (X\otimes \Efun)(V) = |X|\times \Efun(V) \]
for a simplicial set $X$ and $V$ in $\sJ$.
A morphism ``space''
\[ \nat_\pt(\Efun,\Ffun) \]
for arbitrary objects $\Efun$ and $\Ffun$
in $\sE_0$ can be defined as the fibrant simplicial set taking
$[k]$ in $\Delta$ to $\mor_{\sE_0}(\Delta^k\otimes\Efun,\Ffun)$.
This construction does not always respect
weak equivalences in the variables $\Efun$ and $\Ffun$. We therefore
choose a (natural) cofibrant replacement for the source variable, say $\Efun^\natural$ with a
natural weak equivalence $\Efun^\natural\to\Efun$, and put
\[ \Rnat_\pt(\Efun,\Ffun):= \nat_\pt(\Efun^\natural,\Ffun)~. \]
(There is no need to choose a fibrant replacement for $\Ffun$ because $\Ffun$ is already fibrant.)
Then it is easy to verify that $\Rnat_\pt(\Efun,\Ffun)$, viewed as a functor of the first or the
second variable, takes weak equivalences to homotopy equivalences of (fibrant) simplicial sets. \newline
We use the expression \emph{homotopy category} as it is commonly used in connection with model
category structures. For example, a morphism from $\Efun$ to $\Ffun$ in the homotopy category of $\sE_0$
is a homotopy class of morphisms in $\sE_0$ from $\Efun^\natural$ to $\Ffun$. \newline

\begin{exa} \label{exa-humble}
We have $\Rnat_\pt(\Efun,\Ffun)\simeq \Rnat_\pt(T_n\Efun,\Ffun)$
whenever $\Efun,\Ffun$ are objects of $\sE_0$ and
$\Ffun$ is polynomial of degree $\le n$. More precisely, the endofunctor $T_n$ of $\sE_0$ comes with a
natural transformation $\eta_n$ from the identity to $T$. The map
\[ \Rnat_\pt(T_n\Efun,\Ffun) \lra \Rnat_\pt(\Efun,\Ffun) \]
given by composition with $\eta_n\co \Efun\to T_n\Efun$ is a homotopy equivalence.
It is not difficult to produce a proof using further properties of $T_n$ and $\eta_n$ declared in \cite{WeissOrtCal},
but a thorough explanation in model category language can be found in \cite[\S6]{BarnesOman}.
\end{exa}

\medskip
Occasionally we want to view the category of spectra with action of $\Or(n)$ as a model category, too,
in order to make a comparison with the category of homogeneous functors of degree $n$, as a subcategory
of $\sE_0$\,. Most importantly,
the notion of weak equivalence that we use is the coarse one: a morphism between $\Or(n)$-spectra
is a weak equivalence if and only if the underlying morphism of spectra is a weak equivalence.
For the details we can refer to \cite[\S8]{BarnesOman}.

\medskip
Here we are mainly interested in rational phenomena, which simplifies our discussion of
natural transformations. There is a difficulty
that one should be aware of. Let $f\co \Theta\to \Psi$ be a map of spectra with action of $\Or(n)$,
where $n>0$, and let $\bar f\co \Efun_\Theta\to \Efun_\Psi$ be the
corresponding morphism between homogeneous functors of
degree $n$ in $\sE_0$. Thus,
\[ \Efun_\Theta(V)=\Omega^\infty(((V\otimes\RR^n)^c\wedge\Theta)_{\ho\Or(n)}) \]
etc., for all $V$ in $\sJ$. If $f$ is a rational weak equivalence of spectra,
then it does not follow that the map $\bar f\co \Efun_\Theta(V)\to \Efun_\Psi(V)$ is a
rational weak equivalence of spaces for every $V$ in $\sJ$; for example
it may fail to induce a bijection on $\pi_0$. In many cases, however, this is of little
importance to us because of the following principle.

\begin{prop} \label{prop-ratprin}
Let $\Efun$, $\Ffun$ and $\Afun$ be objects of $\sE_0$ which are polynomial of degree $\le n$
for some $n\ge 0$, and let $g\co \Efun\to \Ffun$ be a morphism in $\sE_0$. Suppose that
\begin{itemize}
\item $T_0\Efun(0)\simeq \hocolim_k \Efun(\RR^k)$ is path connected;
\item the map $T_0\Efun(0)\to T_0\Ffun(0)$ induced by $g$ is a rational homotopy equivalence;
\item there exists $\ell\in\ZZ$
such that the derivative spectra $\Theta\Efun^{(i)}$ and $\Theta\Ffun^{(i)}$ are $\ell$-connected
for all $i\ge 1$;
\item $g$ induces a rational weak equivalence $\Theta\Efun^{(i)}\to \Theta\Ffun^{(i)}$ for $1\le i\le n$.
\end{itemize}
Suppose also that $T_0\Afun$ is weakly equivalent to $*$ and
the homotopy groups of the derivative spectra $\Theta^{(i)}\Afun$ are rational vector spaces
throughout. Then $g$ induces a homotopy equivalence
$\Rnat_\pt(\Ffun,\Afun) \lra \Rnat_\pt(\Efun,\Afun)$.
\end{prop}

\proof[Sketch proof] Let $k$ be a non-negative integer. For the purposes of this sketch proof
let $\sJ[k]\subset\sJ$ be the full subcategory spanned by the objects
of dimension at least $k$. A polynomial functor of degree $\le n$ such as $\Efun$ has the property that
the canonical map
\[  \Efun(V) \lra \holimsub{0\ne U\le \RR^{n+1}} \Efun(V\oplus U) \]
is a weak equivalence, for every $V$ in $\sJ$. (Here $U$ runs over the linear subspaces
of $\RR^{n+1}$ and we are using a continuous variant of the homotopy inverse limit.) It follows
that $\Efun$ can be reconstructed, up to weak equivalence, from the restriction of $\Efun$ to $\sJ[1]$,
and by repetition of the argument, from the restriction of $\Efun$ to $\sJ[2]$, or from the restriction
of $\Efun$ to $\sJ[k]$. Therefore
\[ \Rnat_\pt(\Efun,\Afun)\simeq \Rnat_\pt(\Efun_{|\sJ[k]},\Afun_{|\sJ[k]}) \]
holds for arbitrary $k\ge 0$, and similarly for $\Ffun$ in place of $\Efun$. (We are using
a model category structure on the category of continuous functors from $\sJ[k]$ to ${\bf T}_\pt$ which
is analogous to the one that we have on the category of continuous functors from $\sJ$
to ${\bf T}_\pt$\,.) Therefore it is enough to
show that the map
\begin{equation} \label{eqn-bignats}
\Rnat_\pt(\Ffun_{|\sJ[k]},\Afun_{|\sJ[k]}) \lra \Rnat_\pt(\Efun_{|\sJ[k]},\Afun_{|\sJ[k]})
\end{equation}
induced by $g$ is a homotopy equivalence for some $k$, possibly large.
Our conditions on $g$ imply that there exists a natural number $k$ such that
for any $V$ in $\sJ$ of dimension at least $k$, the map
\[  \hofiber[\Efun(V)\to T_0\Efun(V)] \lra \hofiber[\Ffun(V)\to T_0\Ffun(V)] \]
induced by $g$ is a rational weak homotopy equivalence of simply connected based spaces. It follows easily,
by induction on $r$ where $1\le r\le n$, that $g$ induces a homotopy equivalence
\[ \Rmap(\Ffun(V),T_r\Afun(W))\lra \Rmap(\Efun(V),T_r\Afun(W)) \]
when $\dim(V)\ge k$ and $W$ in $\sJ$ is arbitrary. In particular, for $r=n$, we get a homotopy equivalence
\[ \Rmap(\Ffun(V),\Afun(W))\lra \Rmap(\Efun(V),\Afun(W)) \]
induced by $g$, still assuming $\dim(V)\ge k$.
Now we note that there are functors
\[  \Phi_{\Efun,k},\Phi_{\Ffun,k}\co \sJ[k]^\op\times\sJ[k] \lra {\bf T}_\pt \]
given by $(V,W)\mapsto \Rmap(\Efun(V),\Afun(W))$ and $(V,W)\mapsto \Rmap(\Ffun(V),\Afun(W))$,
respectively. We have seen that $g$ induces a weak equivalence $\Phi_{\Ffun,k}\to\Phi_{\Efun,k}$ if
$k$ is sufficiently large. It remains to observe that the spaces in~(\ref{eqn-bignats}) can be
re-defined directly in terms of $\Phi_{\Ffun,k}$ and $\Phi_{\Efun,k}$~, in a homotopy invariant manner,
and that the map (\ref{eqn-bignats}) can be defined in terms of the natural transformation
$\Phi_{\Ffun,k}\to\Phi_{\Efun,k}$ determined by $g$. Such an observation might use a concept like
\emph{homotopy end}. There is a more elementary alternative which relies on explicit
replacements (simplicial resolutions) of $\Efun_{|\sJ[k]}$ and $\Ffun_{|\sJ[k]}$ constructed from the adjoint functor pair
\[
\xymatrix@M=10pt{  {\textup{forget}\co(\textup{cts.~functors from $\sJ[k]$ to ${\bf T}_\pt$})}
\ar@<0.5ex>[r] & \ar@<0.5ex>[l] {({\bf T}_\pt)^{\{k,k+1,k+2,\dots\}}\co \textup{free}}
}
\]
The forgetful functor associates to a functor $\Dfun$ from $\sJ[k]$ to ${\bf T}_\pt$ the sequence of
based spaces $(\Dfun(\RR^i))_{i\ge k}$. The replacements of $\Efun_{|\sJ[k]}$ and $\Ffun_{|\sJ[k]}$ so
obtained are cofibrant replacements
if the values $\Efun(V)$ and $\Ffun(V)$ are already cofibrant in ${\bf T}_\pt$~, for all $V$ in $\sJ[k]$,
which we can assume without loss of generality. \qed

\bigskip
The main point of this section is that our results on the functor $\Sfun_J$ in
proposition~\ref{prop-Rcalc} can be used, together with proposition~\ref{prop-ratprin},
to describe $\Rnat_\pt(\Efun,\Afun)$ in many situations where
$\Efun$ is homogeneous of degree $1$ while $\Afun$ is homogeneous
of degree $2$ and \emph{rational}.

\begin{lem} \label{lem-natfib} Let $\Efun$ and $\Afun$ be objects of $\sE_0$.
Suppose that $\Afun$ is polynomial of degree $\le 2$. Suppose that $T_0\Efun(0)$ is weakly contractible
and that the first and second derivative spectra of $\Efun$ are $\ell$-connected for some $\ell\in\ZZ$.
Then there is a homotopy fiber sequence
\[ \Rnat_\pt(L_2\Efun,\Afun) \leftarrow \Rnat_\pt(\Efun,\Afun)
\leftarrow \Rnat_\pt(T_1\Efun,\Afun).
\]
\end{lem}

\proof By example~\ref{exa-humble} we may replace $\Rnat_\pt(\Efun,\Afun)$
by $\Rnat_\pt(T_2\Efun,\Afun)$. There is a more obvious homotopy fiber sequence
\[ \Rnat_\pt(L_2\Efun,\Afun) \leftarrow \Rnat_\pt(T_2\Efun,\Afun) \leftarrow \Rnat_\pt(\Mfun,\Afun)
\]
where $\Mfun$ is the mapping cone of the canonical morphism $L_2\Efun\to T_2\Efun$. Furthermore we have
$\Rnat_\pt(\Mfun,\Afun) \simeq \Rnat_\pt(T_2\Mfun,\Afun)$
by example~\ref{exa-humble} once again. Therefore it is enough to show that the
natural map $\Mfun\to T_1\Efun$ induces an equivalence of second Taylor approximations,
\[  T_2\Mfun \to T_2T_1\Efun\simeq T_1\Efun~. \]
From the explicit form of the operator $T_2$~, it suffices to show that
$\Mfun(V)\to T_1\Efun(V)$ is $(3d-c)$-connected for all $V$,
where $d=\dim(V)$ and $c$ is a constant independent of $V$. This is a special case of the following
observation related to the Blakers-Massey homotopy
excision theorem: If $f\co Y\to Z$
is a based map where $Z$ is $k$-connected and $f$ is $m$-connected, then the canonical map
from the mapping cone of $\hofiber(f) \to Y$ to $Z$
is $(k+m-c)$ connected. We apply this with $f$ equal to the map $T_2\Efun(V)\to T_1\Efun(V)$
and $k=d-c_1$, $m=2d-c_2$ for suitable constants $c_1$ and $c_2$ which depend on the
integer $\ell$ in our assumptions.
\qed

\medskip
Now suppose that $\Efun=\Sfun_J$ in lemma~\ref{lem-natfib}. Then
two of the spaces in the homotopy fiber sequence of lemma~\ref{lem-natfib} are easy to understand.
The space $\Rnat_\pt(L_2\Sfun_J,\Afun)$ can be understood since we understand natural transformations
between homogeneous functors of the same degree (see lemma~\ref{lem-homfunsamedeg}). The space
$\Rnat_\pt(\Sfun_J,\Afun)$ can be understood because
$\Sfun_J$ behaves in many ways like a (co)representable functor, as shown in the following lemma.

\begin{lem}\label{lem-funcrepre}
Let $\sJ^{|1}$ be the full subcategory of $\sJ$ spanned by the objects $0$ and $\RR$. 
The functor $\Sfun_J$ is freely generated by its restriction  to $\sJ^{|1}$. In particular,
for any $\Ffun$ in $\sE_0$ we have
\[ \Rnat_\pt(\Sfun_J,\Ffun)
\simeq \prod_j\Omega^{k_j}\,\hofiber[\,\Ffun(0)\to \Ffun(\RR).]
\]
\end{lem}

\proof
We show this in the case $\Sfun_J=\Sfun$ (the general case follows similarly).
Let $e$ be the inclusion of $\sJ^{|1}$ in $\sJ$.
By the free generation statement we mean that $\Sfun$ is the left Kan extension of its
restriction  $\Sfun\circ e$. More explicitly, given $V$ in $\sJ$ and a point $y\in \Sfun(V)$
we can find $A$ in $\sJ^{|1}$, an $x\in \Sfun(A)$ and $f\co A\to V$ such that $f_*(x)=y$.
The triple $(A,x,f)$ is unique up to the obvious relations. In the cases where $A=\RR$, we can
always choose $x$ and $f$ such
that  $x\in [0,\infty]$. Then we see that
$\nat_\pt(\Sfun,\Ffun)$ can be identified with the space of pairs of based maps $(f,g)$ making the diagram
\[
 \xymatrix{
\{0,\infty\} \ar[r]^-{\textup{incl}} \ar[d]_f & [0,\infty] \ar[d]^g \\
\Ffun(0) \ar[r]^-{\textup{incl}_*} & \Ffun(\RR)
}
\]
commute, where $\{0,\infty\}=\Sfun(0)$ and $[0,\infty]\subset \Sfun(\RR)$.
This space of pairs $(f,g)$ is just the homotopy fiber of the inclusion-induced map $\Ffun(0)\to \Ffun(\RR)$.
\newline
Furthermore, it is easy to see from this description that $\Sfun$ is cofibrant
in the specified model category structure on $\sE_0$.
Therefore our computation of $\nat_\pt(\Sfun,\Ffun)$ can be taken as a computation of $\Rnat_\pt(\Sfun,\Ffun)$.
\qed

\begin{lem}\label{lem-homfunsamedeg}
Let $\Efun$ and $\Ffun$ be two homogeneous functors in $\sE_0$\,, of the same degree $n>0$.
Then $\Rnat_\pt(\Efun,\Ffun)\simeq (\Rmap(\Theta \Efun^{(n)},\Theta \Ffun^{(n)}))^{\ho\Or(n)}$.
\end{lem}

\proof This is suggested by the classification of homogeneous functors~\cite{WeissOrtCal},
but for a really thorough proof see \cite{BarnesOman}. \qed

\begin{lem}\label{lem-homogmorph}
Suppose that $\Afun$ in $\sE_0$ is the homogeneous functor of degree $2$ corresponding
to an Eilenberg-MacLane spectrum $\underline{H}\QQ$ with trivial action of $\Or(2)$.
Let $\Efun$ in $\sE_0$ be any homogeneous functor of degree $1$ such that
the first derivative spectrum of $\Efun$ is
$0$-connected, and each of its homotopy groups is finitely generated.
Then $\Rnat_\pt(\Efun,\Afun)$ is contractible.
\end{lem}

\proof
We start the proof with a special case.
Suppose that $\Efun=T_1\Sfun_J$ where $k_j>0$ for all $j=1,2,...$~.
By lemma~\ref{lem-funcrepre} we have
\[
\begin{array}{ccl}
\Rnat_\pt(\Sfun_J,\Afun)
&\simeq & \prod_{j}\Omega^{k_j}\hofiber[\Afun(0)\to \Afun(\RR)] \\
&\simeq & \prod_{j}\Omega^{k_j}\Omega^\infty\hofiber[(S^0\wedge \underline{H}\QQ)_{\ho\Or(2)}\to
(S^2\wedge \underline{H}\QQ)_{\ho\Or(2)}] \\
&\simeq & \prod_{j}\Omega^{k_j+1}\Omega^\infty\left(((S^2/S^0)\wedge \underline{H}\QQ)_{\ho\Or(2)}\right) \\
&\simeq & \prod_{j}\Omega^{k_j}\Omega^\infty\left((S^1_+\wedge \underline{H}\QQ)_{\ho\Or(2)}\right) \\
&\simeq & \prod_{j}\Omega^{k_j}\Omega^\infty\left((S^0\wedge \underline{H}\QQ)_{\ho\Or(1)}\right) \\
&\simeq & \prod_{j}\Omega^{k_j}\QQ
\end{array}
\]
which is contractible. (We have used: $S^2/S^0$ as a based $\Or(2)$-space
is homeomorphic to a smash product $S^1\wedge S^1_+$ with $\Or(2)$
acting trivially on the first factor $S^1$, and by the standard nontrivial action on the second
factor $S^1_+$~.) By lemma~\ref{lem-homfunsamedeg}  we have
\[\Rnat_\pt(L_2\Sfun_J,\Afun)\simeq \map^{\ho\Or(2)}(\Theta \Sfun_J^{(2)},\underline{H}\QQ)\]
and by proposition~\ref{prop-Rcalc} the right-hand side is contractible. From the homotopy fiber sequence
of lemma~\ref{lem-natfib} it follows that $\Rnat_\pt(T_1\Sfun_J,\Afun)$ is contractible. \newline
Now for the general case: Let $\Efun'$ be the homogeneous functor of degree 1 corresponding
to the rationalization $\Theta_\QQ$ of the first derivative spectrum $\Theta$ of $\Efun$. Then
\[ \Rnat_\pt(\Efun,\Afun)\simeq \Rnat_\pt(\Efun',\Afun) \]
by inspection or by proposition~\ref{prop-ratprin}. But $\Theta_\QQ$
is a retract, in the homotopy category of spectra with action of $\Or(1)$, of $\Psi_\QQ$
where $\Psi$ is the first derivative spectrum of $\Sfun_J$,
for some sequence $J$ as in the previous step. Let $\Efun''$ be the homogeneous functor
of degree $1$ associated with $\Psi_\QQ$. Then $\Rnat_\pt(\Efun',\Afun)$ is
a retract up to homotopy of $\Rnat_\pt(\Efun'',\Afun)$ and the latter is
homotopy equivalent to $\Rnat_\pt(T_1\Sfun_J,\Afun)$ by proposition~\ref{prop-ratprin}.
But $\Rnat_\pt(T_1\Sfun_J,\Afun)$ is contractible as we have seen. \qed

\medskip
The next lemma is a generalization of lemma~\ref{lem-homogmorph} using very similar ideas.

\begin{lem}\label{lem-homogmorph2}
Let $\Afun$ be as in lemma~\ref{lem-homogmorph}. Let
$\Efun$ and $\Ffun$ in $\sE_0$ be homogeneous functors of degree $1$.
Suppose that the first derivative spectra of $\Efun$ and $\Ffun$
are $0$-connected and $(-1)$-connected, respectively, with finitely generated homotopy groups.
Then the map from $\Rnat_\pt(\Efun\times\Ffun,\Afun)$ to $\Rnat_\pt(\Ffun,\Afun)$ induced
by the inclusion $\Ffun\to \Efun\times\Ffun$ is a homotopy equivalence.
\end{lem}

\proof As in the proof of lemma~\ref{lem-homogmorph}, we can quickly reduce to the case where
$\Efun$ is $T_1\Sfun_J$ and $\Ffun$ is $T_1\Sfun_H$ for suitable sequences $J$ and $H$. Then
$\Efun\times\Ffun\simeq T_1(\Sfun_J\vee \Sfun_H)$. By the homotopy fiber sequence of
lemma~\ref{lem-natfib}, it is enough to show that the following restriction maps (both
induced by the inclusion $\Sfun_H\to \Sfun_J\vee\Sfun_H$) are
homotopy equivalences:
\[  \Rnat_\pt(L_2(\Sfun_J\vee\Sfun_H),\Afun) \lra \Rnat_\pt(L_2\Sfun_H,\Afun) \]
\[  \Rnat_\pt(\Sfun_J\vee\Sfun_H,\Afun) \lra  \Rnat_\pt(\Sfun_H,\Afun)~.\]
For the first one of these we use proposition~\ref{prop-Rcalc} and lemma~\ref{lem-homfunsamedeg}.
For the second, we note
\[ \Rnat_\pt(\Sfun_J\vee\Sfun_H,\Afun)\simeq \Rnat_\pt(\Sfun_J,\Afun)\times \Rnat_\pt(\Sfun_H,\Afun) \]
where $\Rnat_\pt(\Sfun_J,\Afun)$ is contractible, as seen in the proof of lemma~\ref{lem-homogmorph}. \qed

\medskip
With a view to applications in section~\ref{sec-pessimist}, we ask whether some of
lemma~\ref{lem-homogmorph} remains intact if we allow more choice for the target functor $\Afun$.
An answer is given in lemma~\ref{lem-shifty} below. This relies on the following, which is well known.

\begin{lem} \label{lem-specshift} Let $\Efun$ belong to $\sE_0$~, with $n$-derivative spectrum
$\Theta=\Theta\Efun^{(n)}$. The $n$-th derivative spectrum of $\Efun(-\oplus\RR^k)$ can be identified with
$S^{nk}\wedge\Theta$
where $S^{nk}=(\RR^k\otimes\RR^n)^c$, with the standard action of $\Or(n)$ on $\RR^n$
and the trivial action on $\RR^k$. \emph{(We use the diagonal action of $\Or(n)$ on $S^{nk}\wedge\Theta$.)}
\end{lem}

\proof It is easy to reduce the statement to the case where $\Efun$ is homogeneous of degree $n$.
(If $\Efun$ is not homogeneous of degree $n$, replace it by the homogeneous layer $L_n\Efun$.)
In the case where $\Efun$ is homogeneous of degree $n$, it is weakly equivalent to the functor
\[ V~\mapsto~ \Omega^\infty(((V\otimes\RR^n)^c\wedge\Theta)_{\ho\Or(n)}). \]
Then $\Efun(-\oplus\RR^k)$ is weakly equivalent to the functor
\[ V~\mapsto~ \Omega^\infty(((V\otimes\RR^n)^c\wedge(\RR^k\otimes\RR^n)^c\wedge \Theta)_{\ho\Or(n)})~. \]
Therefore $\Efun(-\oplus\RR^k)$ is also homogeneous of degree $n$ and its derivative spectrum is
$S^{nk}\wedge \Theta$ as claimed. \qed

\begin{lem} \label{lem-shifty}
Suppose that $\Afun$ in $\sE_0$ is the homogeneous functor of degree $2$ corresponding
to an Eilenberg-MacLane spectrum $S^r\wedge\underline{H}\QQ$ with trivial action of $\Or(2)$,
where $r\ge 0$.
Let $\Efun$ in $\sE_0$ be any homogeneous functor of degree $1$ such that
the first derivative spectrum $\Theta$ of $\Efun$ is
$(-\ell)$-connected, where $\ell\ge 0$, and each homotopy group $\pi_k\Theta$ is finitely generated.
Then the map
\[  \Rnat_\pt(\Efun,\Afun) \to \Rnat_\pt(\Efun,\Afun(-\oplus\RR^{\ell+2})) \]
induced by $\Afun\to \Afun(-\oplus\RR^{\ell+2})$ is nullhomotopic.
\end{lem}

\proof We begin with the special case where $\ell=0$. So
$\Theta$ is $0$-connected. As in the proof of lemma~\ref{lem-homogmorph}, we can easily reduce to the case where
$\Efun$ is $T_1\Sfun_J$ for some sequence $J=(k_j)_{j=1,2,\dots}$ of strictly positive integers.
Write $\Afun^\sigma$ for $\Afun(-\oplus\RR)$ and $\Afun^{\sigma\sigma}$ for $\Afun(-\oplus\RR^2)$.
For use later on, note that the second derivative spectrum of $\Afun^\sigma$ is $S^2\wedge S^r\wedge\underline{H}\QQ$\,,
with $\Or(2)$ acting through the standard representation on $S^2=\RR^2\cup\infty$\,; similarly
the second derivative spectrum of $\Afun^{\sigma\sigma}$ is $S^2\wedge S^2\wedge S^r\wedge\underline{H}\QQ$,
with $\Or(2)$ acting diagonally through the standard representation on both copies of $S^2=\RR^2\cup\infty$.
See lemma~\ref{lem-specshift}.
We have a commutative diagram
\[
\xymatrix@R=15pt{
\Rnat_\pt(L_2\Sfun_J,\Afun) \ar[d]&\ar[l] \Rnat_\pt(T_2\Sfun_J,\Afun) \ar[d]^-u &\ar[l]\ar[d]
\Rnat_\pt(T_1\Sfun_J,\Afun) \\
\Rnat_\pt(L_2\Sfun_J,\Afun^\sigma) \ar[d]^-v  & \ar[l] \Rnat_\pt(T_2\Sfun_J,\Afun^\sigma)
\ar[d] & \ar[d]\ar[l] \Rnat_\pt(T_1\Sfun_J,\Afun^\sigma) \\
\Rnat_\pt(L_2\Sfun_J,\Afun^{\sigma\sigma}) & \ar[l]
\Rnat_\pt(T_2\Sfun_J,\Afun^{\sigma\sigma}) & \ar[l] \Rnat_\pt(T_1\Sfun_J,\Afun^{\sigma\sigma})
}
\]
where the rows are homotopy fiber sequences by lemma~\ref{lem-natfib}. The vertical arrows are induced
by $\Afun\to \Afun^\sigma$ and $\Afun^\sigma\to \Afun^{\sigma\sigma}$. Therefore, if we can produce a
nullhomotopy for the arrow labelled $u$, and another nullhomotopy for the arrow labeled $v$,
then we have a nullhomotopy for the composition in the right-hand column. \newline
In order to make a nullhomotopy for $u$, we calculate as in the proof of lemma~\ref{lem-homogmorph}:
\[  \begin{array}{rcl}
\Rnat_\pt(T_2\Sfun_J,\Afun) &\simeq\cdots\simeq & \prod_j \Omega^{k_j}(S^r\wedge \underline{H}\QQ) \\
\Rnat_\pt(T_2\Sfun_J,\Afun^\sigma) & \simeq\cdots\simeq & \prod_j \Omega^{k_j}(S^2\wedge S^r\wedge\underline{H}\QQ)~.
\end{array}
\]
The map $u$ corresponds to the inclusion
\[ S^r\wedge\underline{H}\QQ~\cong~S^0\wedge S^r\wedge\underline{H}\QQ \lra S^2\wedge S^r\wedge\underline{H}\QQ \]
and so we have a nullhomotopy for it. \newline
To make a nullhomotopy for $v$, we write $\Psi$ for the second derivative spectrum of $\Sfun_J$ and
use lemma~\ref{lem-homfunsamedeg} first to rewrite $v$ as a map
\[
\xymatrix@R=12pt{
(\Rmap(\Psi,S^2\wedge S^r\wedge \underline{H}\QQ))^{\ho\Or(2)} \ar[d] \\
(\Rmap(\Psi,S^2\wedge S^2\wedge S^r\wedge \underline{H}\QQ))^{\ho\Or(2)}.
}
\]
This is induced by the inclusion
\[ S^2\wedge S^r\wedge \underline{H}\QQ~\cong~S^0\wedge S^2\wedge S^r\wedge \underline{H}\QQ
~\lra~S^2\wedge S^2\wedge S^r\wedge \underline{H}\QQ~.\]
It follows from proposition~\ref{prop-Rcalc} that
$\Psi$ is weakly equivalent to a spectrum, with action of $\Or(2)$, of the form
\[  \Or(2)_+\wedge_{\Sigma_2}\Phi \]
where $\Phi$ is a spectrum  with an action of the symmetric group $\Sigma_2$\,. Therefore
\[
\begin{array}{rcl}
(\Rmap(\Psi,S^2\wedge S^r\wedge \underline{H}\QQ))^{\ho\Or(2)} & \simeq &
(\Rmap(\Phi,S^2\wedge S^r\wedge \underline{H}\QQ))^{\ho\Sigma_2} \\
(\Rmap(\Psi,S^2\wedge S^2\wedge S^r\wedge\underline{H}\QQ))^{\ho\Or(2)} & \simeq &
(\Rmap(\Phi,S^2\wedge S^2\wedge S^r\wedge\underline{H}\QQ))^{\ho\Sigma_2}~.
\end{array}
\]
Now it is clear that $v$ admits a nullhomotopy. \newline
In the general case, $\ell\ge 0$,
we can use the following commutative diagram:
\[
\xymatrix{
\Rnat_\pt(\Efun,\Afun) \ar[r] \ar[d] &  \Rnat_\pt(\Efun,\Afun(-\oplus\RR^{\ell+2})) \\
\Rnat_\pt(\Efun(-\oplus\RR^\ell),\Afun(-\oplus\RR^\ell)) \ar[r] &
\Rnat_\pt(\Efun(-\oplus\RR^\ell),\Afun(-\oplus\RR^{\ell+2})) \ar[u]
}
\]
The first derivative spectrum of $\Efun(-\oplus\RR^\ell)$ is $S^\ell\wedge\Theta$, which is
$0$-connected. By what we have already shown, the lower horizontal arrow admits a nullhomotopy. \qed

\medskip
The next lemma is about spaces of unbased natural transformations. For such a space we write
$\Rnat(\Cfun,\Dfun)$, where $\Cfun$ and $\Dfun$ can be objects of $\sE_0$ or of $\sE$,
the category of continuous functors from $\sJ$
to $\bf T$. We regard $\sE$ as a model category by analogy with $\sE_0$\,, so that the meaning
of $\Rnat(\Cfun,\Dfun)$ should be clear.

\begin{lem} \label{lem-homogmorph3} In the situation of lemma~\ref{lem-homogmorph2}, the
inclusion $\Ffun\to \Efun\times\Ffun$ induces a homotopy
equivalence $\Rnat(\Efun\times\Ffun,\Afun)\lra\Rnat(\Ffun,\Afun)$.
\end{lem}

\proof There is a commutative diagram
\[
\xymatrix{
\Rnat_\pt(\Efun\times\Ffun,\Afun) \ar[r] \ar[d] & \Rnat_\pt(\Ffun,\Afun) \ar[d] \\
\Rnat(\Efun\times\Ffun,\Afun) \ar[r] \ar[d]_{\textup{eval. at base pt}} &
\Rnat(\Ffun,\Afun) \ar[d]^{\textup{eval. at base pt}} \\
\Afun(0) \ar@{=}[r] & \Afun(0)
}
\]
in which the columns are homotopy fiber sequences. It follows (with lemma~\ref{lem-homogmorph2})
that the inclusion $\Ffun\to \Efun\times\Ffun$ induces a homotopy
equivalence from the portion of $\Rnat(\Efun\times\Ffun,\Afun)$ lying above the base element of $\pi_0\Afun(0)$
to the portion of $\Rnat(\Ffun,\Afun)$ lying above the base element of $\pi_0\Afun(0)$. Let $y\in\Afun(0)$ be a point
representing some other element of $\pi_0\Afun(0)$, if possible.
The functor $\Afun$ with the new base point $y\in \Afun(0)$ is a new object of $\sE_0$ which is still
homogeneous of degree $2$ and still has second derivative spectrum $\underline{H}\QQ$, with trivial action of
$\Or(2)$. (The concept of homogeneous functor has a base-point free definition. The $m$-th derivative spectra of any functor
$\Dfun$ in $\sE$ form a fibered family of spectra on the space $T_0\Dfun(0)=\hocolim_i\Dfun(\RR^i)$; if
$\Dfun$ is homogeneous of degree $m>0$, then $T_0\Dfun(0)$ is weakly contractible.) Therefore we
may conclude that $\Ffun\to \Efun\times\Ffun$ induces a homotopy
equivalence from the portion of $\Rnat(\Efun\times\Ffun,\Afun)$ lying above any particular
element of $\pi_0\Afun(0)$
to the portion of $\Rnat(\Ffun,\Afun)$ lying above that same element of $\pi_0\Afun(0)$. \qed

\begin{cor}\label{cor-homogmorph3}
With $\Afun$ as in lemma~\ref{lem-homogmorph}, suppose that $g\co \Dfun\to\Ffun$
is a morphism in $\sE_0$ where $\Ffun$ and $\Dfun$ are polynomial of degree $\le 1$
and $T_0g\co T_0\Dfun(0)\to T_0\Ffun(0)$ is a weak equivalence of \emph{path connected} spaces. Suppose
that the first derivative spectra of $\Dfun$ and $\Ffun$ are $(-1)$-connected and that
the map between them induced by $g$
is 1-connected, and admits a right inverse in the homotopy category of spectra with action of
$\Or(1)$. Then $g$ induces a homotopy equivalence
$\Rnat(\Ffun,\Afun) \lra \Rnat(\Dfun,\Afun)$.
\end{cor}

\proof We view the statement as a parameterized version of lemma~\ref{lem-homogmorph3}.
The path connected space $Y:=T_0\Ffun(0)\simeq T_0\Dfun(0)$ serves as the parameter space.
Let $\simp(Y)$ be the category of singular simplices
of $Y$. An object is a pair $(m,\sigma)$ where $\sigma\co\Delta^m\to Y$ is a singular $m$-simplex.
A morphism from $(m,\sigma)$  to $(n,\psi)$ is a monotone map $u\co [m]\to[n]$ such that $\sigma=\psi\circ u_*$~.
For $y=(m_y,\sigma_y)$ in $\simp(Y)$ let $\Ffun_y$ be the homotopy pullback of
\[
\xymatrix{
\Delta^m \ar[r]^\sigma & Y \ar[r] & T_0\Ffun & \ar[l] \Ffun
}
\]
and let $\Dfun_y$
be the homotopy pullback of
\[
\xymatrix{
\Delta^m \ar[r]^\sigma & Y \ar[r] & T_0\Ffun & \ar[l]  \Ffun   & \ar[l]_g \Dfun
}
\]
(where we have taken the liberty to view spaces as constant functors on $\sJ$).
Since $Y$ can be identified (up to a weak equivalence) with the homotopy colimit of
the functor $y\mapsto \Delta^{m_y}$ on $\simp(Y)$, it follows that $\Ffun$ can be identified
with the homotopy colimit (in $\sE$) of $y\mapsto \Ffun_y$
and $\Dfun$ can be identified with the homotopy colimit of $y\mapsto \Dfun_y$\,.
These homotopy colimits are turned into homotopy inverse limits
by $\Rnat(-,\Afun)$. Therefore the horizontal arrows in the following commutative diagram
\[
\xymatrix{
\Rnat(\Ffun,\Afun) \ar[d]^{g^*} \ar[r] & \holim_y~\Rnat(\Ffun_y,\Afun) \ar[d]^{g^*} \\
\Rnat(\Dfun,\Afun) \ar[r] & \holim_y~\Rnat(\Dfun_y,\Afun)
}
\]
are homotopy equivalences. Therefore, in order to show that the left-hand vertical arrow is a
homotopy equivalence, we only need to show that the right-hand vertical arrow is a homotopy
equivalence. We can achieve that by showing that
\[ \Rnat(\Ffun_y,\Afun) \lra \Rnat(\Dfun_y,\Afun) \]
induced by $g$ is a homotopy equivalence for every $y$ in $\simp(y)$. Since $Y$ is path
connected, it is enough to show this when $y$ is the base point of $Y$. In that case
$\Ffun_y=L_1F$ and $\Dfun_y=L_1\Dfun$. We get what we need from lemma~\ref{lem-homogmorph2}
and our assumption that $L_1\Dfun\to L_1\Ffun$ admits a splitting. \qed

\section{Splitting hypotheses} \label{sec-split}
\begin{prop}\label{prop-hypc}
Hypothesis~\thehypc is equivalent to saying that the inclusion-induced
$\Or(2)$-map $\Theta \Ofun^{(2)}\to \Theta \Tfun^{(2)}$ admits a
rational weak left inverse.
\end{prop}

(By a rational weak left inverse we mean an $\Or(2)$-map $\,\Theta\Tfun^{(2)}\to \Theta'$ of spectra
with an action of $\Or(2)$ such that the composite map
$\Theta \Ofun^{(2)}\to \Theta \Tfun^{(2)}\to \Theta'$
is a rational weak homotopy equivalence.) \newline
One half of proposition~\ref{prop-hypc} is trivial. If $\Tfun\to\Efun$ is a rational
weak left inverse (over $B\ZZ/2$) for
the inclusion $\Ofun\to\Tfun$, then the map of second derivative spectra
induced by $\Tfun\to \Efun$ is a rational weak left inverse for the map of second derivative spectra induced by $\Ofun\to\Tfun$.
The proof of the other half occupies the remainder of this section.

\medskip
We recall how a functor $\Dfun$ in $\sE_0$ determines a sequence of spectra
$\Theta\Dfun^{(i)}$,
for $i=1,2,3,\dots$ (for more details on this and what follows, see~\cite{WeissOrtCal}
and \cite{BarnesOman}). The category $\sJ$ is contained in a larger category $\sJ_i$ enriched over based spaces and
the functor $\Dfun$ has a right Kan extension 
\[ \Dfun^{(i)}\co \sJ_i \to {\bf T}_\pt\,. \] 
An explicit formula for $\Dfun^{(i)}$ is
\begin{equation} \label{eqn-rightKan}  \Dfun^{(i)}(V)= \hofiber~\bigl[
\,\,\Dfun(V) \lra \!\holimsub{\quad 0\ne U\le \RR^i} \Dfun(V\oplus U)\,\,\bigr]
\end{equation}
where we use a topologically enhanced homotopy limit.
Instead of saying that
$\Dfun^{(i)}$ is defined on $\sJ_i$ we can also pretend that it is defined on $\sJ$ and comes with the
following additional structure: a natural map
\begin{equation}\label{eqn-susp}  \sigma\co (V\otimes\RR^i)^c\wedge \Dfun^{(i)}(W) \lra \Dfun^{(i)}(V\oplus W)
\end{equation}
with the expected associativity and unital properties. (Sometimes it is more convenient to write
$\Dfun^{(i)}(W\oplus V)$ for the target; the two options are canonically identified.)
Moreover $\Dfun^{(i)}$ comes with an action of $\Or(i)$, obvious
from the explicit formula, such
that $\sigma$ is equivariant. (It is equivariant for the diagonal action of $\Or(i)$ on the source. By specializing to $V=\RR$
and $W=\RR^i$ we obtain a spectrum with twisted action of $\Or(i)$ where the structure maps involve smash product
with a sphere $(\RR^i)^c$ on which $\Or(i)$ acts nontrivially. This can be untwisted.
Besides, it is essential in the following that we don't
specialize too soon.)

\begin{rem} \label{rem-towersplit}
We need to state and explain a few facts to be used later.
\begin{itemize}
\item[(a)]  The inclusion $B\Or\to B\TOP$ is a rational homotopy equivalence.
\item[(b)]  The canonical actions of the group $\Or$ on $\Theta \Ofun^{(i)}$ and $\Theta \Tfun^{(i)}$ for $i\ge1$
are special cases of a  natural action of $\Or$ on all spectra. (That natural action of $\Or$ on all
spectra extends to a well-known natural action of $G\subset QS^0$ on all spectra.)
\end{itemize}
Statement (a) was mentioned in the introduction. It can be restated as
saying that the inclusion of zeroth Taylor approximations $T_0\Ofun \to T_0\Tfun$ is a rational homotopy equivalence.\newline
Regarding (b), suppose that $\Dfun$ is a continuous functor from
$\sJ$ to based spaces and that $\Dfun(\RR^\infty)=\hocolim_n \Dfun(\RR^n)$ is path-connected.
Then the $i$-th derivative spectrum of $\Dfun$ is defined (for $i=1,2,3,\dots$) and
comes with an action of $\Or(i)$. What matters here is that it also comes with an action of $\Omega \Dfun(\RR^{\infty})$,
commuting with the action of $\Or(i)$ (modulo replacement of the $i$-th derivative spectrum by something weakly equivalent
in the category of spectra with action of $\Or(i)$). Equivalently, the $i$-th derivative spectrum
of $\Dfun$ is just one fiber of a fibered spectrum over $\Dfun(\RR^\infty)$, with fiberwise action of $\Or(i)$.
This is part of the general theory \cite{WeissOrtCal}, but we need to recall how it works. Fix $n\ge 0$ and $x\in \Dfun(\RR^n)$.
Let $s_n\Dfun_x$ be defined by $s_n\Dfun_x(V)=\Dfun(\RR^n\oplus V)$, with base point equal to the image of
$x$ under the inclusion-induced map $\Dfun(\RR^n)\to \Dfun(\RR^n\oplus V)$. The functor $s_n\Dfun_x$ in $\sE_0$
determines derivative spectra
\begin{equation} \label{eqn-fibspectra}
\Theta(s_n\Dfun_x)^{(i)}
\end{equation}
for $i\ge 1$. As $x$ runs through $\Dfun(\RR^n)$ these spectra constitute a fibered spectrum
over $s_n\Dfun(0)=\Dfun(\RR^n)$, with
fiberwise action of $\Or(i)$. Now let us see how this fibered spectrum over $\Dfun(\RR^n)$
depends on $n$. Let $x\in \Dfun(\RR^n)$ and let $y\in \Dfun(\RR^{n+1})$ be the image
of $x$ under the map induced by the inclusion $\RR^n\to \RR^{n+1}$.
Then, using the maps~(\ref{eqn-susp}) with $V=\RR$, we get a homotopy equivalence
\[  \Theta(s_n\Dfun_x)^{(i)} \lra \Omega^{\RR^i}\Theta(s_{n+1}\Dfun_y)^{(i)} \]
which we can also write in the form
\begin{equation} \label{eqn-fibspec}
  \Omega^{\RR^n\otimes\RR^i}\Theta(s_n\Dfun_x)^{(i)} \lra \Omega^{\RR^{n+1}\otimes\RR^i}\Theta(s_{n+1}\Dfun)^{(i)}_y~.
  \end{equation}
Summarizing, the family $x\mapsto \Omega^{ni}\Theta(s_n\Dfun_x)^{(i)}$ where $x\in \Dfun(\RR^n)$ is a fibered spectrum
over $\Dfun(\RR^n)$ and the maps~(\ref{eqn-fibspec}) allow us to assemble these fibered spectra, by a telescope construction,
to a (quasi-)fibered spectrum over $\Dfun(\RR^\infty)=\hocolim_n \Dfun(\RR^n)$.
This quasi-fibered spectrum over
$\Dfun(\RR^\infty)$, for $i=1,\dots,k$,
is one of the ingredients in a stagewise description of the $k$-th Taylor approximation $T_k\Dfun$ of $\Dfun$.
That is how we will use it below, with $k=2$.

We now specialize by taking $\Dfun=\Ofun$, while $i$ remains unspecified. 
In order to reason carefully we use a specific model of $\Ofun$. Write $Q=\RR^\infty=\colim_n\RR^n$.
Let $\Ofun(V)$ be the space of linear subspaces of $V\oplus Q$ of the same dimension as $V$ (topologized as a colimit
of ordinary finite-dimensional and compact Grassmannians).
This is a continuous functor of the variable $V$ in $\sJ$. Key observation: for $U\in \Ofun(\RR^n)$,
there is a zigzag of equivalences relating the functors $s_n\Ofun_U$ and $s_U\Ofun_U$ in $\sE_0$\,, where $s_U\Ofun_U(W):=\Ofun(U\oplus W)$
with the standard base point. Indeed, for $W$ of dimension $k$ we have based homotopy equivalences
\[
\CD
{s_n\Ofun_U(W)=\{\textup{$(n+k)$-dimensional linear subspaces of $\RR^n\oplus W\oplus Q$}\}} \\
@VV {\textup{embed $Q$ as first copy of two}} V \\
{\{\textup{$(n+k)$-dimensional linear subspaces of $\RR^n\oplus W\oplus Q\oplus Q$}\}} \\
@AA {\textup{embed $Q$ as second copy of two}} A \\
{s_U\Ofun_U(W)=\{\textup{$(n+k)$-dimensional linear subspaces of $U\oplus W\oplus Q$}\}}
\endCD
\]
(In the top line, the base point is the sum of $W$ and the copy of $U$ in $\RR^n\oplus 0\oplus Q$; in the bottom line,
it is simply $U\oplus W=U\oplus W\oplus 0$ as a linear subspace of $U\oplus W\oplus Q$.)
Therefore
\[ \Omega^{ni}\Theta(s_n\Ofun_U)^{(i)}\simeq\cdots\simeq \Omega^{ni}\Theta(s_U\Ofun_U)^{(i)}\simeq
\Omega^{ni}((U\otimes\RR^i)^c\wedge\Theta\Ofun^{(i)}). \]
This should be seen as a zigzag of equivalences of fibered spectra over $\Ofun(\RR^n)$. Therefore we have identified the monodromy
of this fibered spectrum (an action of $\Omega\Ofun(\RR^n)$ on the fiber over the base point) as the action of $\Or(n)$ on
\[ \Omega^{ni}((\RR^n\otimes\RR^i)^c\wedge\Theta\Ofun^{(i)})  \]
induced by the standard action of $\Or(n)$ on $\RR^n$, hence on $\RR^n\otimes\RR^i$ and on the one-point
compactification of $\RR^n\otimes\RR^i$. We now let $n$ tend
to infinity (the details of that limit process are left to the reader)
and thereby complete our sketch proof of (b) in the case of $\Ofun$. \newline
The case of $\Tfun$ is similar provided we make a good start. What is required
most of all is an explicit model of $\Tfun$ sufficiently similar to the one for $\Ofun$ which we have just seen.
It is not easy to give a description of $\Tfun(V)$ as a kind of
topological Grassmannian. We prefer to describe $\Tfun(V)$ as the classifying space of a topological
groupoid which we call $A(V)$ for now. The object space of the groupoid $A(V)$ is $\Ofun(V)$, that is, the space of linear subspaces of $V\oplus Q$ of
the same dimension as $V$. The morphism space of $A(V)$ is the space of triples $(W_1,W_2,h)$ where $W_1,W_2$ are linear subspaces of $V\oplus Q$ of
the same dimension as $V$ and $h\co W_1\to W_2$ is a homeomorphism. The source of a morphism $(W_1,W_2,h)$ is $W_1$ and the target
is $W_2$. The nerve of $A(V)$ is a simplicial \emph{space} and the classifying space $BA(V)$
is the geometric realization of that. Evidently $A(V)$ is equivalent to a topological group, and therefore the
natural inclusion of the true classifying space of the topological group $\TOP(V)$ into $BA(V)$ is a weak equivalence.
Last not least, this explicit model of $\Tfun$ lends itself to a proof of (b) for the spectra $\Theta\Tfun^{(i)}$ which is analogous
to the above proof of (b) for the spectra $\Theta\Ofun^{(i)}$. Note once again that we are not interested in the action of
$\TOP\simeq \Omega \Tfun(\RR^\infty)$ on the $i$-th derivative spectrum
of $\Tfun$, but only in the action of the subgroup $\Or \simeq \Omega \Ofun(\RR^\infty)$.
\end{rem}

\begin{rem} \label{rem-TaylorPost} The Taylor tower in orthogonal calculus (of a functor $\Dfun$ from $\sJ$
to ${\bf T}_\pt$) has strong formal
similarities with the Postnikov tower of a based, connected CW-space $X$. The essentially constant
functor $T_0\Dfun$ plays the part of $B\pi_1(X)$ in Postnikov theory. The $i$-th derivative spectrum of $\Dfun$
plays a role similar to that of the homotopy group $\pi_{i+1}=\pi_{i+1}(X)$. The fact that $\pi_{i+1}$ is a module over
$\pi_1(X)$ is analogous to the fact that $\Theta^{(i)}$ extends to a fibered spectrum
\[ z \mapsto \Theta^{(i)}_z \]
on the space $T_0\Dfun(0)=\Dfun(\RR^\infty)$.
The inductive construction of the stages $X_{i+1}$ of the Postnikov tower of $X$ is best described
by means of homotopy pullback squares
\[
\xymatrix{
X_{i+1} \ar[rr]^-{\textup{proj.}} \ar[d]^-{\textup{proj.}} && X_i \ar[d]^-{ \kappa_i}\\
B\pi_1(X) \ar[rr]^-{\textup{zero section}} &&(B^{i+2}\pi_{i+1})_{\ho \pi_1(X)}~.
}
\]
Here the lower right-hand term is the homotopy orbit construction for the action of $\pi_1(X)$
on $B^{i+2}\pi_{i+1}$, so that there is a projection from
it to $B\pi_1(X)$ with Eilenberg-MacLane fiber $B^{i+2}\pi_{i+1}$. By analogy with that, there is a
homotopy pullback square
\[
\xymatrix{
T_i\Dfun \ar[rr]^-{\textup{proj.}} \ar[d]^-{\textup{proj.}} && T_{i-1}\Dfun \ar[d]^-{ \kappa_i}\\
T_0\Dfun \ar[rr]^-{\textup{zero section}} && \Hfun{[S^1\wedge\Theta^{(i)}_\bullet]}~.
}
\]
Here $\Hfun{[S^1\wedge\Theta^{(i)}_\bullet]}$ is a functor which is essentially
determined by the space $T_0\Dfun(0)$ and
the fibered spectrum
\[ z\mapsto S^1\wedge \Theta^{(i)}_z \]
on it, with the fiberwise action of $\Or(i)$. There is a
forgetful projection
\[ \Hfun{[S^1\wedge\Theta^{(i)}_\bullet]} \lra T_0\Dfun\,, \]
and for every point $z\in T_0\Dfun(0)$ the homotopy fiber of that projection at $z$ is the homogeneous
functor of degree $i$ from $\sJ$ to ${\bf T}_\pt$ determined by the spectrum
$S^1\wedge \Theta^{(i)}_z$
with action of $\Or(i)$.
\end{rem}

\proof[Proof of prop.~{\rm\ref{prop-hypc}}]
To begin we replace $\Tfun$ by a functor $\TOfun$, rationally equivalent over $B\ZZ/2$
to $\Tfun$. This is defined as the
homotopy pullback of
\[
\xymatrix@R=12pt{
& \Tfun \ar[d] \\
T_0\Ofun \ar[r] & T_0\Tfun
}
\]
By proposition~\ref{prop-calcE} it is enough to produce a rational splitting
(rational left inverse in a suitable homotopy category) for the inclusion
\[  T_2\Ofun\lra T_2\TOfun~. \]
We rely on remark~\ref{rem-TaylorPost} and so imagine a commutative diagram with two horizontal
arrows to be constructed:
\[
\xymatrix@R=12pt{
T_2\Ofun \ar[d] \ar[r] & T_2\TOfun \ar[d] \ar@{..>}[r]^v & \Efun \ar[d] \\
T_1\Ofun \ar[d] \ar[r] & T_1\TOfun \ar[d] \ar@{..>}[r]^u & T_1\Efun \ar[d] \\
T_0\Ofun \ar[r]^-\simeq & T_0\TOfun \ar[r]^-\simeq & T_0\Efun
}
\]
The composite map in the lower row and the composite map in the middle row are meant to be weak equivalences.
The composite map in the top row is meant to be a rational weak equivalence. \newline
We begin with the construction of the arrow $u$. From the point of view of remark~\ref{rem-TaylorPost},
the functors
$T_1\Ofun$ and $T_1\TOfun$ are determined by the space
\[ T_0\Ofun(0)=T_0\TOfun(0)=B\Or \]
and fibered spectra
$\Theta^1_\bullet$\,, $\Psi^1_\bullet$ over $B\Or$ (the families of first derivative spectra of $\Ofun$ and $\TOfun$),
and sections (generic name $\kappa_1$) of the fibrations
\[
\begin{array}{lcr}
\Omega^\infty((S^1\wedge\Theta^1_\bullet)_{\ho\Or(1)}) &\lra& B\Or \\
\Omega^\infty((S^1\wedge\Psi^1_\bullet)_{\ho\Or(1)}) &\lra& B\Or
\end{array}
\]
respectively. Therefore it is enough to show (for the construction of arrow $u$)
that the map (over $B\Or$) of fibered first derivative spectra $\Theta^1_{\bullet}\to\Psi^1_{\bullet}$
induced by $\Ofun\to \TOfun$ admits a homotopy left inverse. This is clear from proposition~\ref{prop-calcF}
and part (b) of remark~\ref{rem-towersplit}. We also use the splitting, with $\Or(1)$-invariance,
of $\underline{A}(*)$ into the sphere spectrum and another wedge summand. This can be obtained by looking
at the first derivative spectra in $\Ofun\to\Tfun\to\Gfun$. See proposition~\ref{prop-calcxi}. \newline
We come to the construction of the arrow $v$.
From the point of view of remark~\ref{rem-TaylorPost},
the functor $T_2\Ofun$ is determined by the functor $T_1\Ofun$\,, a fibered spectrum $\Theta^2_{\bullet}$ on $B\Or$
with fiberwise action of $\Or(2)$ and a natural transformation
\[ T_1\Ofun \lra \Hfun[S^1\wedge\Theta^2_{\bullet}] \]
(generic name $\kappa_2$)
over $T_0\Ofun$. Similarly, $T_2\TOfun$ is determined by $T_1\TOfun$\,, a fibered spectrum $\Psi^2_{\bullet}$ on $B\Or$
with fiberwise action of $\Or(2)$ and a natural transformation
\[ T_1\TOfun \lra \Hfun[S^1\wedge\Psi^2_{\bullet}] \]
(generic name $\kappa_2$) over $T_0\TOfun\simeq T_0\Ofun$. (To avoid confusion, beware that $\Hfun[...]$ is
a construction which depends on the following input: \emph{an integer $i>0$} and a fiberwise spectrum on a space,
with fiberwise action of $\Or(i)$. The integer $i$ is not shown in our notation; here $i=2$.)
By our assumption and part (b) of
remark~\ref{rem-towersplit}, there exist a fibered spectrum $\Xi_\bullet$ over $B\Or$ with fiberwise
action of $\Or(2)$, and a map of fibered spectra $p\co \Psi^2_\bullet\to \Xi_\bullet$ preserving the fiberwise
actions of $\Or(2)$, and serving as a rational left homotopy inverse
for the map $\Theta^2_\bullet\to \Psi^2_\bullet$ induced by the inclusion $\Ofun\to \TOfun$. We may as well
assume that the fibers of $\Xi_\bullet$ are spectra whose homotopy groups are rational vector spaces;
then it follows that they are all homotopy equivalent to $\Omega\underline{H}\QQ$. More precisely,
the fibered spectrum $\Xi_\bullet$ is fiber homotopy trivial,
i.e., each $\Xi_y$ can be identified with $\Omega\underline{H}\QQ$,
with the trivial action of $\Or(2)$. As a result $\Hfun[S^1\wedge\Xi_{\bullet}]$ is weakly equivalent to
a product $T_0\TOfun\times\Afun$, where  $\Afun$ is the homogeneous functor of degree $2$
associated with the spectrum $\underline{H}\QQ$ and the trivial action of $\Or(2)$ on it.
The composition
\[
\xymatrix{
T_1\TOfun \ar[r]^-{\kappa_2} & \Hfun[S^1\wedge\Psi^2_{\bullet}] \ar[r]^p & \Hfun[S^1\wedge\Xi_{\bullet}]
}
\]
(which is a morphism over $T_0\TOfun$) then amounts to a morphism $T_1\TOfun\to \Afun$ in $\sE_0$\,.
And therefore the construction of arrow $v$ amounts to finding a factorization (in the homotopy category of $\sE_0$)
of that morphism $T_1\TOfun\to \Afun$ in the following way:
\[
\xymatrix{
T_1\TOfun \ar[r]^u & T_1\Efun \ar@{..>}[r] & \Afun~.
}
\]
Such a factorization exists by corollary~\ref{cor-homogmorph3} (set $g=u$ in~\ref{cor-homogmorph3}).  \qed

\section{Orthogonal calculus and smoothing theory}\label{sec-orcalsmooth}
The main input from smoothing theory is the following general theorem due to Morlet \cite{Morlet}.
See also the thorough exposition of Morlet's result in \cite{KirbySiebenmann} and the earlier
work \cite{HirschMazur}, which was later to appear in print.

\begin{thm} \label{thm-Morlet}
The space of smooth structures on a closed topological manifold $M$ of dimension $m\ne 4$ is
homotopy equivalent, by an obvious forgetful map, to the space of vector bundle structures on the topological
tangent (micro-)bundle of $M$. \newline
For a compact topological $m$-manifold $M$ with smooth boundary, $m\ne 4$,
the space of smooth structures on $M$ extending the given structure on $\partial M$ is
homotopy equivalent to the space of vector bundle structures on the topological
tangent bundle of $M$ extending the prescribed vector bundle structure over $\partial M$.
\end{thm}

\begin{rem}
There is a homotopy lifting principle for vector bundle structures
on fibre bundles with fiber $\RR^m$. Namely, if $E\to X\times [0,1]$
is such a fiber bundle and a vector bundle structure has been chosen
on $E|_{X\times 0}$,  then this vector bundle structure admits an
extension to all of $E$. This has the following homotopy theoretic
consequence. Given a map $c:X\to B\TOP(m)$, the space of vector
bundle structures on the associated fiber bundle on $X$ with fiber
$\RR^m$ is homotopy equivalent to the space of maps $\tilde c$ from
$X\times [0,1]$ to $B\TOP(m)$ which satisfy $\tilde c(x,0)=c(x)$ and
map $X\times 1$ to $B\Or(m)\subset B\TOP(m)$.
\end{rem}

\begin{exa}\label{exa-reghofi} Let $\sY_n$ be the space of smooth structures on $D^n$ extending the standard
smooth structure on $S^{n-1}$. Then
\[ \sY_n~\simeq~\Omega^n(\TOP(n)/\Or(n))~. \]
Furthermore, there is a homotopy fiber sequence
\[ \reg \lra \Omega^2\sY_n \lra \sY_{n+2} \]
where $\reg=\reg(n,2)$ as previously defined. It is obtained as
follows: Any smooth regular map $f\co D^n\times D^2\to D^2$
satisfying our boundary conditions is a smooth fiber bundle by
Ehresmann's theorem. The underlying bundle of topological manifolds
is canonically trivial relative to the given trivialization over the
boundary $\partial D^2$. (Its structure group, the group of
topological automorphisms of $D^n$ extending the identity on the
boundary, is contractible by the Alexander trick.) Hence $f$
determines a family, parametrized by $D^2$, of smooth structures on
$D^n$ extending the standard smooth structure on $S^{n-1}$. This
family is of course trivialized over the boundary $\partial D^2$.
The resulting ``integrated'' smooth structure on the total space of
the bundle is equal to the standard structure on $D^n\times D^2$ by
assumption. These observations lead to the stated homotopy fiber
sequence, and we conclude
\begin{equation} \label{eqn-modelreg}
\begin{split}
\reg~\simeq~\Omega^{n+2}\hofiber[\TOP(n)/\Or(n)\to\TOP(n+2)/\Or(n+2)]~.
\end{split}
\end{equation}
\end{exa}

\begin{exa}\label{exa-smoothhcobdisc} Let $\mathscr{H}_n$ be the
space of smooth structures on $D^{n-1}\times [0,1]$ extending the standard structure on
$(D^{n-1}\times 0)\,\cup\,( \partial D^{n-1}\times [0,1])$. Reasoning as in the previous example we
have  a homotopy equivalence
\[ \mathscr{H}_{n}~\simeq~\Omega^{n-1}\big(\hofiber\big[\TOP(n-1)/\Or(n-1)\rightarrow\TOP(n)/\Or(n)\big]\big)~, \]
if $n\geq 6$.
The space $\mathscr{H}_n$ can also be viewed as the space of smooth h-cobordisms on $D^{n-1}$, since the space of
topological h-cobordisms on $D^{n-1}$ is contractible by the Alexander trick and the affirmed Poincar\'e conjecture.
\end{exa}

\begin{exa} With the notation of the previous examples, there is a homotopy fiber sequence
\[ \aut_{\diff}(D^n) \lra \aut_{\topo}(D^n) \lra \sY_n \]
where $\aut_{\diff}(D^n)$ is the space of diffeomorphisms $D^n\to D^n$ which
extend the identity of $S^{n-1}$, and $\aut_{\topo}(D^n)$ is the
topological analogue. This homotopy fiber sequence is obtained by considering the action
of $\aut_{\topo}(D^n)$ on $\sY_n$, and the stabilizer subgroup of the base point in $\sY_n$.
By the Alexander trick, $\aut_{\topo}(D^n)$ is contractible. Hence
$\aut_{\diff}(D^n)~\simeq~\Omega^{n+1}(\TOP(n)/\Or(n))$.
\end{exa}

We turn to the relation between the \thehypb and \thehypc hypotheses. Smoothing theory
allows us to reformulate hypothesis \thehypb so that it fits into the orthogonal calculus framework.

\begin{rem}\label{rem-untwist}
Let $\Gamma$ be a compact Lie group. We are interested in based spaces and spectra with an action of the group $\Gamma$.
As in section~\ref{sec-homoext}, it is convenient to assume that the based spaces are objects of ${\bf T}_\pt$~.
Therefore, let ${\bf T}_\pt^G$ be the category of such spaces with an action of $G$.
We use the model category structure on ${\bf T}_\pt^G$ where a morphism $X\to Y$ is a weak equivalence
(or a fibration) if the underlying morphism in ${\bf T}_\pt$ is a weak equivalence (or a fibration). \newline
Many examples of spectra with an action of $\Gamma$ which
arise in orthogonal calculus have some features which are reminiscent of an equivariant setting,
with notions of stability involving possibly nontrivial representations $V$ of $G$.
We strive to suppress such features. The following questions arise frequently:
\begin{enumerate}
\item[(i)] Given a fixed representation $V$ of $\Gamma$ (by which we mean an object $V$ of $\sJ$
and an action of $\Gamma$ on it), a sequence of based $\Gamma$-spaces $(X_{nV})_{n\in\NN}$ and based $\Gamma$-maps
from $V^c\wedge X_{nV}$ to $X_{(n+1)V}$ (with the diagonal action on the source), can we build a
spectrum $\underline{X}$ with action of $\Gamma$ from these data? Also,
given two such sequences $(X_{nV})_{n\in\NN}$ and $(Y_{nV})_{n\in\NN}$, and compatible based $\Gamma$-maps
$f_n\co X_{nV}\to Y_{nV}$, can we build a $\Gamma$-map $f\co \underline{X}\to \underline{Y}$~?
\item[(ii)] In these circumstances, can the $\Gamma$-map $f$ be recovered if we only know the based $\Gamma$-maps
$\Omega^{nW}f_n\co \Omega^{nW}X_{nV} \to \Omega^{nW}Y_{nV}$ for all $n$, where $W$ is another representation of $\Gamma$~?
We are willing to assume that the $Y_{nV}$ are Eilenberg-MacLane spaces and the adjoints of the maps
$V^c\wedge Y_{nV}\to Y_{(n+1)V}$ are homotopy equivalences from $Y_{nV}$ to $\Omega^VY_{(n+1)V}$.
\end{enumerate}
The following propositions try to answer these questions.
\end{rem}

\begin{prop}\label{prop-untwistGpectra}
Given a fixed representation $V$ of $\Gamma$, a sequence of based $\Gamma$-spaces $(X_{nV})_{n\in\NN}$ and based $\Gamma$-maps
$V^c\wedge X_{nV}\to X_{(n+1)V}$, the spaces
\[ \underline{X}(m) := \hocolimsub{n\to\infty}\Omega^{nV}(S^m\wedge X_{nV}) \]
(for $m\ge 0$) form an $\Omega$-spectrum $\underline{X}$. \emph{(More precisely, maps
$\underline{X}(m)\to \Omega\underline{X}(m)$ will be defined which are weak homotopy equivalences.)}
\end{prop}

\proof There are obvious structure maps
\begin{equation} \label{eqn-untwistinf} \hocolimsub{n\to\infty}\Omega^{nV}(S^m\wedge X_{nV}) \to
\Omega\bigl(\hocolimsub{n\to\infty}\Omega^{nV}(S^{m+1}\wedge X_{nV})\bigr)
\end{equation}
We need to show that these are weak homotopy equivalences. Suppose to begin with that
there exists $n_0$ such that, for all $n\ge n_0$, the $\Gamma$-map
$V^c\wedge X_{nV}\to X_{(n+1)V}$ is a homeomorphism. Then
$X_{nV}$ is $(\dim(nV)-k)$-connected for a constant $k$ independent of $n$. It follows by Freudenthal's theorem
that the maps~(\ref{eqn-untwistinf}) are weak homotopy equivalences, since they can be written in the form
\[
\hocolimsub{n\to\infty}\Omega^{nV}(S^m\wedge X_{nV}) \to
\hocolimsub{n\to\infty}\Omega^{nV}\left(\Omega(S^{m+1}\wedge X_{nV})\right).
\]
The general case follows from this special case by a direct limit argument. \qed

\medskip
Let $Y$ be a based $\Gamma$-space (where $\Gamma$ is a compact Lie group as before).
We denote by $Y_{\ho\Gamma}$ the (unreduced) Borel construction as usual. However in the case of
a spectrum $\underline{X}$ with action of $\Gamma$ as in proposition~\ref{prop-untwistGpectra},
we write $\underline{X}_{\ho\Gamma}$ for the spectrum obtained by applying
the reduced Borel construction levelwise.

\begin{prop}\label{prop-untwistreconst}
Keeping the notation of proposition~\ref{prop-untwistGpectra}, let $W$ be another
representation of $\Gamma$ and assume
$\Gamma$ is connected. If $d=\dim(V)$ is greater than $e=\dim(W)$, the canonical homomorphisms
\[ g^n\co H^{k}(\underline{X}_{\ho \Gamma};\QQ)\to  H^{k+nd-ne}_{\Gamma}
(\Omega^{nW}X_{nV},\pt;\QQ)~,\]
induce an injection
\[ g\co  H^{k}(\underline{X}_{\ho \Gamma};\QQ)\to \prod_n H^{k+nd-ne}_{\Gamma}
(\Omega^{nW}X_{nV},\pt;\QQ)~.\]
\end{prop}

\proof By the universal coefficient theorem in rational homology,
it is enough to show that the corresponding homology homomorphisms $g_n$ produce a surjection
\[ \bigoplus_n H_{k+nd-ne}^\Gamma(\Omega^{nW}X_{nV},\pt;\QQ)
\lra H_{k}(\underline{X}_{\ho \Gamma};\QQ)~. \]
We note that $\oplus_ng_n$ is defined  by the composition
\begin{equation}\label{eqn-cohospec}
\begin{split}
\xymatrix@R=15pt{
\bigoplus_n H_{k+nd-ne}^\Gamma(\Omega^{nW}X_{nV},\pt;\QQ) \ar[d] \\
\bigoplus_n H_{k+nd}^\Gamma(S^{nW}\wedge\Omega^{nW}X_{nV},\pt;\QQ)\ar[d]\\
{\colimsub{n} H_{k+nd}^\Gamma(X_{nV},\pt;\QQ) \cong  H_{k}(\underline{X}_{\ho \Gamma};\QQ)}~,
}
\end{split}
\end{equation}
where the top arrow is the suspension isomorphism and the bottom one is induced by the obvious maps
$S^{nW}\wedge\Omega^{nW}X_{nV}\to X_{nV}$. (Instead of $(\RR^n\otimes W)^c$ or $(nW)^c$
we have written $S^{nW}$.)
Now we need to show the bottom arrow is onto. By a direct limit argument (as in the previous proposition) we can reduce this
to the case where $X_{nV}$ is $(nd-p)$-connected for a constant $p$ independent of $n$. It then follows that the map
\[ S^{nW}\wedge\Omega^{nW}X_{nV}\to X_{nV}\]
is $(2nd-ne-q)$-connected for some constant $q$ independent of $n$.
For large enough $n$ we have
$2nd-ne-q= n(d-e)-q+nd>k+nd$,
so that the $n$-th summand in the middle term of diagram~(\ref{eqn-cohospec}) maps onto
$H_{k+nd}^\Gamma(X_{nV},\pt;\QQ)$.
\qed

\begin{rem} In proposition~\ref{prop-untwistreconst}, we can interpret an element $f$ of
$H^k(\underline{X}_{\ho \Gamma};\QQ)$ as a homotopy class of
$\Gamma$-maps from a cofibrant replacement of $\underline{X}$ to $\underline{Y}\simeq S^k\wedge H\QQ$.
We may also assume that
$\underline{Y}$ is constructed from based $\Gamma$-spaces $Y_{nV}$ and $\Gamma$-maps
\[ V^c\wedge Y_{nV}\to Y_{(n+1)V} \]
whose adjoints are homotopy equivalences $Y_{nV}\to \Omega^VY_{(n+1)V}$. The image of $f$ in the group
$H^{k+nd-ne}_\Gamma(\Omega^{nW}X_{nV},\pt;\QQ)$
is the map $\Omega^{nW}f_n$ from $\Omega^{nW}X_{nV}$ to $\Omega^{nW}Y_{nV}$ induced by $f$.
The content of the proposition is that $f$ is
determined by these images $\Omega^{nW}f_n$~. This gives an affirmative answer to the question in (ii)
of remark~\ref{rem-untwist}.
\end{rem}

\begin{rem}
The proof of proposition~\ref{prop-untwistreconst} proves more than what is stated.
In fact, given any infinite subset $S$ of the
natural numbers, the composition of the injection $g$ with the projection
\[ \prod_{n\in \NN} H^{k+nd-ne}_\Gamma(\Omega^{nW}X_{nV},\pt;\QQ)
\lra \prod_{n\in S} H^{k+nd-ne}_\Gamma(\Omega^{nW}X_{nV},\pt;\QQ)
\]
is still an injection.
\end{rem}

\begin{prop} \label{prop-evenfromhyptohyp} Hypothesis \thehypb implies
that the map of second derivative spectra induced by the inclusion $\Ofun\to\Tfun$
admits a rational weak left inverse. \emph{(Compare proposition~\ref{prop-hypc}.)}
\end{prop}

\proof Let $\Theta_{\Ofun}=\Theta \Ofun^{(2)}$, $\Theta_{\Tfun}=\Theta \Tfun^{(2)}$ and
$\Theta_{\Ofun\to \Tfun}=\hofiber[\Theta_{\Ofun}\to \Theta_{\Tfun}]$.
We need to show that the forgetful map $\Theta_{\Ofun\to \Tfun}\to\Theta_{\Ofun}$ is rationally
nullhomotopic with $\Or(2)$-invariance (after cofibrant replacement inflicted on the source).
As $\Theta_{\Ofun}$ is rationally an Eilenberg-MacLane spectrum concentrated in dimension $-1$ with
trivial $\Or(2)$-action,
this amounts to showing that a class $\delta$ in
the spectrum cohomology $H^{-1}((\Theta_{\Ofun\to \Tfun})_{\ho \Or(2)};\QQ)$ vanishes.
A transfer argument shows that the homomorphism
\[ H^{-1}((\Theta_{\Ofun\to \Tfun})_{\ho \Or(2)};\QQ)\lra H^{-1}((\Theta_{\Ofun\to \Tfun})_{\ho S^1};\QQ)\]
determined by restriction is injective. Therefore we only have to show that $\delta$ is zero in
the group $H^{-1}((\Theta_{\Ofun\to \Tfun})_{\ho S^1};\QQ)$, or equivalently
that the forgetful map $\Theta_{\Ofun\to \Tfun}\to\Theta_{\Ofun}$ is rationally
nullhomotopic with weak $S^1$-invariance.
This is equivalent to showing that the forgetful map
$\Omega^V\Theta_{\Ofun\to \Tfun}\to\Omega^V\Theta_{\Ofun}$
is rationally nullhomotopic with weak $S^1$-invariance, where $V$ is the standard
2-dimensional representation of $S^1$.
Again this is equivalent to the vanishing of a cohomology class
$\eta\in H^{-3}((\Omega^V\Theta_{\Ofun\to \Tfun})_{\ho S^1};\QQ)$.
By proposition~\ref{prop-untwistreconst} this will follow from showing that the cohomology classes
$\eta_n\in H^{-3+2n-n}_{S^1}(\Omega^{n}X_{nV},\pt;\QQ)$
determined by $\eta$ are zero for all even $n$, where
\[ X_{nV}=\Omega^V\hofiber[\Ofun^{(2)}(\RR^n)\to \Tfun^{(2)}(\RR^n)]. \]
To show this we use the commutative diagram of $S^1$-spaces
\begin{equation} \label{eqn-lastgasp}
\begin{split}
\xymatrix@C=45pt@R=15pt{ \Omega^n\Omega^V\hofiber[\Ofun^{(2)}(\RR^n)\to \Tfun^{(2)}(\RR^n)] \ar[r]^-{\eta_n}\ar[d] &
\Omega^n\Omega^V \Ofun^{(2)}(\RR^n) \ar[d]^-\tau \\
 \Omega^n\Omega^V\hofiber\left[\frac{\Or(n+2)}{\Or(n)}\to \frac{\TOP(n+2)}{\TOP(n)}\right] \ar[r]\ar[d]^{\simeq} &
{\Omega^n\Omega^V\left(\frac{\Or(n+2)}{\Or(n)}\right)_\QQ} \\
\reg(n,2) \ar[r]^{\nabla} & {\Omega^n\Omega^V\left(\frac{\Or(n+2)}{\Or(n)}\right)_\QQ} \ar@{=}[u]
}
\end{split}
\end{equation}
The map $\tau$ is a rational homotopy equivalence.
(Indeed, there is a homotopy fiber sequence
\[ \CD \Ofun^{(2)}(\RR^n) @>\tau >>  \displaystyle\frac{\Or(n+2)}{\Or(n)} @>>> \Gamma(E\to\RR P^1)~, \endCD \]
where $E\to\RR P^1$ is the fiber bundle with fiber given by $E_L=\Or(\RR^{n+2})/\Or(\RR^n\oplus L)$
and \linebreak $\Gamma(E\to\RR P^1)$ is the corresponding section space.
It is easy to check that the base point component of $\Omega^n\Omega^V\Gamma(E\to\RR P^1)$ is
rationally contractible.)
Therefore by hypothesis~\ref{hypb}, all these classes $\eta_n$ are zero.
\qed

\begin{rem} Our computations in section~\ref{sec-ortho} show that the map of second derivative
spectra induced by $\Ofun\to \Tfun$ does admit a rational weak left inverse as a map
of spectra (\emph{action of $\Or(2)$ suppressed}).
Indeed, by proposition~\ref{prop-calcxi} the map of second derivative
spectra induced by the composition $\Ofun\to \Tfun\to \Gfun$ admits a rational weak left inverse as a map
of spectra (action of $\Or(2)$ suppressed).
\end{rem}

\section{Pessimistic option} \label{sec-pessimist} Here we
look for weaker versions of hypotheses \thehypa, \thehypb and \thehypc.
Weaker versions might hold if the original
hypotheses turn out to be wrong. Or we can be ambitious in pessimism, trying to disprove even
the weaker versions.

\begin{thm} \label{thm-pessimist} The following are equivalent.
\begin{itemize}
\item[a)] There exists an even positive integer $k$ such that
for all even $n\ge 0$, the Pontryagin class in $H^{2n+2k}(B\STOP(n);\QQ)$ is decomposable.
\item[b)] There exists an even positive $k$ such that for all even $n\ge 4$, the class
\[ c^k\cup[\nabla] \in H_{S^1}^{n-3+2k}(\reg(n,2),\pt;\QQ) \]
is zero, where $c\in H^2(BS^1;\ZZ)$ is the standard generator.
\item[c)] There exists an even positive $k$ such that
the inclusion
$\Theta\Ofun^{(2)}\lra S^{2k}\wedge \Theta\Ofun^{(2)}$
of spectra with action of $\Or(2)$ admits a rational weak factorization
through the map $\Theta\Ofun^{(2)}\to \Theta\Tfun^{(2)}$ induced by the
inclusion $\Ofun\to \Tfun$.
\end{itemize}
\end{thm}

\smallskip
\emph{Comment on }a). Since $B\STOP(n)$ is simply connected, saying that the
Pontryagin class in $H^{2n+2k}(B\STOP(n);\QQ)$ is decomposable amounts to saying that it evaluates to zero on any element
of $\pi_{2n+2k}B\STOP(n)$. It is also equivalent to saying that the looped class in the group
$H^{2n+2k-1}(\STOP(n);\QQ)$ is zero.

\smallskip
\emph{Comment on }b). For any pair of spaces $(X,Y)$ with action of a Lie group $G$, the Borel cohomology
$H^*_G(X,Y;\QQ)$
is a graded module over the graded ring $H^*(BG;\QQ)$. This uses the projections
from $X_{\ho G}$ and $Y_{\ho G}$ to $BG$.

\smallskip
\emph{Comment on }c). Think of $S^{2k}$ as $(\RR^k\otimes\RR^2)^c$ where $\Or(2)$ acts on the
factor $\RR^2$ by the standard action. The fixed point set of this action of $\Or(2)$ on $S^{2k}$
is identified with $S^0$. Use the diagonal action of $\Or(2)$ on
$S^{2k}\wedge \Theta\Ofun^{(2)}$.
The inclusion of the fixed point set, $S^0\to S^{2k}$, induces
\[ \Theta\Ofun^{(2)}~\cong~S^0\wedge\Theta\Ofun^{(2)}\lra S^{2k}\wedge \Theta\Ofun^{(2)}. \]
By \emph{rational weak factorization} etc.~we mean a
commutative square of spectra with action of $\Or(2)$
\begin{equation} \label{eqn-weakfact2}
\begin{split}
\xymatrix@C=35pt@R=20pt{
\Theta\Ofun^{(2)} \ar[r]^-{\textup{incl.}} \ar[d]^-{\textup{incl.}} & S^{2k}\wedge\Theta\Ofun^{(2)}
\ar@{..>}[d]^{\simeq_\QQ}  \\
\Theta\Tfun^{(2)} \ar@{..>}[r] & \underline{E}
}
\end{split}
\end{equation}
where the arrow labeled $~\simeq_\QQ~$ is a rational weak equivalence (of spectra).

\smallskip
\emph{Comment on the integer $k$.} It is not claimed that a), b) and c) are equivalent for the same
fixed $k$. See remark~\ref{rem-funnyk} below.

\medskip
As a warm-up for showing a)$\Rightarrow$b) we make some remarks on the meaning of a).
It follows from a) that there is a rational factorization up to homotopy of the following kind:
\begin{equation} \label{eqn-warm1}
\begin{split}
\
\xymatrix@u@C=40pt@R=10pt{
\Omega^n(\Or/\Or(n)) \ar[r]^-{\Omega^{n+1}p_{(n+k)/2}} \ar[d]  & K(\QQ,n+2k-1) \\
\Omega^n(\TOP/\TOP(n)) \ar@{..>}[ur]
}
\end{split}
\end{equation}
In saying that, we have not mentioned actions of $\SOr(2)=S^1$. But it does not
cost us anything to view~(\ref{eqn-warm1}) as a diagram of based $S^1$-spaces. We write
\[
\begin{array}{rcl}
\Or/\Or(n) & = & \colim_N~\Or(\RR^n\oplus\RR^2\oplus\RR^N)/\Or(\RR^n), \\
\TOP/\TOP(n) & = & \colim_N~\TOP(\RR^n\oplus\RR^2\oplus\RR^N)/\TOP(\RR^n) \\
\end{array}
\]
and let $S^1$ act via conjugation (by rotations) on the summand $\RR^2$. Also, $S^1$ acts trivially
on $K(\QQ,n+2k-1)$. The $\Omega^n$ prefix is attached afterwards, with
a trivial action of $\Or(2)$ on the loop coordinates. It does not cost us anything because these actions
are trivial up to weak equivalence, as we noted earlier in the proof of proposition~\ref{prop-AtoB}.

\proof[Proof of $\textup{a})\Rightarrow\textup{b})$] This is very similar to the proof of
\thehypa$\Rightarrow$\thehypb in proposition~\ref{prop-AtoB}. We use the notation introduced there,
but we assume a) instead of \thehypa, with a fixed even $k>0$ and $n\ge 4$, and we are going to deduce b) with the same
$k$ and $n$\,. We write $Y_n$ and $Z_n$ instead of $Y$ and $Z$ to be more specific; $Y_n$ means
$\GL(n+2)/\GL(n)$ or preferably $\Or(n+2)/\Or(n)$, and $Z_n$ means $\TOP(n+2)/\TOP(n)$.
In addition we will need
\[  W_n=\Ofun^{(2)}(\RR^n) = \hofiber[\,\Ofun(\RR^n)\to \holimsub{0\ne U\le \RR^2}\Ofun(\RR^n\oplus U)\,]. \]
(We use a topological homotopy inverse limit.) This is from the sequence of spaces which makes up the
second derivative spectrum of $\Ofun$ in its original untwisted form. There is a forgetful map $W_n\to Y_n$. \newline
Since the composition~(\ref{eqn-lose}) represents zero in the homotopy category of based $S^1$-spaces,
it is enough to show that there exists a class
\[  y \in H^{n-3+2k}_{S^1}(\Omega^{n+2}Z_n,\pt;\QQ) \]
which under the inclusion
$\Omega^{n+2}Y_n\lra \Omega^{n+2}Z_n$
pulls back to a cup product
\[ c^k\cup\gamma\in H^{n-3+2k}_{S^1}(\Omega^{n+2}Y_n,\pt;\QQ) \]
for some nonzero element
$\gamma \in H^{n-3}_{S^1}(\Omega^{n+2}(Y_n)_\QQ,\pt;\QQ)~\cong~\QQ$.
To that end we set up a commutative diagram of based $S^1$-spaces
\begin{equation} \label{eqn-pessimist1}
\begin{split}
\xymatrix@u@C=25pt@R=-45pt{
{\Omega^{n+2}W_n} \ar[r] \ar[dr]_-{\textup{forget}}  &  {\Omega^{n+2}W_{n+k}}  \ar[dr]^{\textup{forget}}  \\
& {\Omega^{n+2}Y_n\rule{4mm}{0mm}} \ar[r] \ar[dr]_-{\textup{inc.}} \ar[dd]_-{\textup{inclusion}} &
{\Omega^{n+2}Y_{n+k}\rule{4mm}{0mm}} \ar[dr]^-{\textup{inc.}} \\
  & & {\rule{3mm}{0mm}\Omega^{n+2}(\Or/\Or(n))} \ar[r] \ar[dd] & {\Omega^{n+2}(\Or/\Or(n+k))}
\ar[rr]^-{\Omega^{n+3}p_{(n+k)/2}} && K(\QQ,n+2k-3)  \\
& \Omega^{n+2}Z_n \ar[dr]_-{\textup{inc.}}  \\
&  & \Omega^{n+2}(\TOP/\TOP(n)) \ar@{..>}[rrruu]
}
\end{split}
\end{equation}
The dotted arrow comes from diagram~(\ref{eqn-warm1}), and strictly speaking it makes that triangle or quadrilateral
commutative in a suitable (rational) homotopy category of $S^1$-spaces only.  Note that we have inflicted $\Omega^2$
with the standard (nontrivial) action of $S^1=\SOr(2)$ on the two loop coordinates. \newline
Since $W_n$ and $W_{n+k}$ are constituents of the second derivative spectrum of the functor $\Ofun$
in the original twisted form, there is a
structure map
\begin{equation} \label{eqn-pessimist2}
S^{2k}\wedge W_n \lra W_{n+k}~.
\end{equation}
It is an $S^1$-map,
where $S^1$ acts on $S^{2k}=(\RR^k\otimes\RR^2)^c$ via its standard action by rotations on $\RR^2$,
and diagonally on the entire source. By proposition~\ref{prop-calcE}, the second derivative spectrum
of the functor $\Ofun$ is rationally an $\Omega$-spectrum, that is, $W_n\simeq_{\QQ} K(\QQ,2n-1)$
and $W_{n+k}\simeq_\QQ K(\QQ,2n+2k-1)$ and the map
$W_n\to \Omega^{2k}W_{n+k}$
adjoint to~(\ref{eqn-pessimist2}) is a rational homotopy equivalence. The $S^1$-map
\[ \Omega^{n+2}W_n\to \Omega^{n+2}W_{n+k} \]
in diagram~(\ref{eqn-pessimist1}) is simply the restriction of~(\ref{eqn-pessimist2})
to $W_n~\cong~S^0\wedge W_n$, with $\Omega^{n+2}$ inflicted afterwards.
Therefore, under that map $\Omega^{n+2}W_n\to \Omega^{n+2}W_{n+k}$, the fundamental cohomology class
\[  \Omega^{n+2}v_{n+k} \in H^{n+2k-3}_{S^1}(\Omega^{n+2}W_{n+k},\pt;\QQ) \]
pulls back to the cup product $c^k\,\cup\,\Omega^{n+2}v_n\in H^{n+2k-3}_{S^1}(\Omega^{n+2}W_n,\pt;\QQ)$,
up to multiplication by a nonzero scalar. On the other hand we know (proof of proposition~\ref{prop-evenfromhyptohyp}, discussion of the map $\tau$)
that the forgetful arrow from $W_n$ to $Y_n$
in (\ref{eqn-pessimist1}) induces a rational homotopy equivalence
$\Omega^{n+2}W_n \to \Omega^{n+2}(Y_n)_\QQ$.
This means that we can solve our problem by a diagram chase in diagram~(\ref{eqn-pessimist1}).
We start with the fundamental cohomology class of the Eilenberg-MacLane space at the
top of the diagram (but we regard it as a class in Borel cohomology $H^*_{S^1}$)
and pull it back to $\Omega^{n+2}Z_n$ to obtain $y$ with the required
properties. \qed

\proof[Proof of \textup{b)$\Rightarrow$c)}]
Our assumption b) is equivalent to the statement
that the composition of $S^1$-maps
\[
\xymatrix{
 \reg(n,2) \ar[r] & \Omega^{n+2} Y_\QQ \lra S^{2k}\wedge \Omega^{n+2}Y_\QQ
}
\]
is nullhomotopic as a based $S^1$-map, for all even $n\ge 4$, after cofibrant replacement of the source
$\reg(n,2)$. Here $Y=\GL(n+2)/\GL(n)$ and $S^{2k}$ is $(\RR^k\otimes\RR^2)^c$
with the standard action of $\Or(2)$
on $\RR^2$. As in the proof of proposition~\ref{prop-evenfromhyptohyp},
and with the same abbreviations, this implies that the composition of $\Or(2)$-maps
$\Theta_{\Ofun\to \Tfun}\to\Theta_{\Ofun} \hookrightarrow S^{2k}\wedge \Theta_{\Ofun}$
is rationally nullhomotopic as an $\Or(2)$-map (after cofibrant replacement inflicted on the source).
That in turn is equivalent to the statement that the inclusion
$\Theta_{\Ofun} \hookrightarrow S^{2k}\wedge \Theta_{\Ofun}$,
as an $\Or(2)$-map, admits a weak (rational) factorization through the map
from $\Theta_{\Ofun}$  to $\Theta_{\Tfun}$
determined by the inclusion $\Ofun\to \Tfun$. \qed

\smallskip
\proof[Proof of $\textup{c)}\Rightarrow\textup{a)}$] Let $L_{[1,2]}\Ofun$ and $L_{[1,2]}\Tfun$
be the functors defined by
\[  \hofiber[ T_2\Ofun\to T_0\Ofun]~,\qquad
\hofiber[T_2\Tfun\to T_0\Tfun] \]
respectively.
We are going to show that c) for a particular even
$k$ implies a commutative diagram in the homotopy category of $\sE_0$~,
\[
\xymatrix{
L_{[1,2]}\Ofun \ar[r] \ar[d] &  L_{[1,2]}\Ofun(-\oplus\RR^\ell) \ar@{..>}[d]_{\simeq_\QQ} \\
L_{[1,2]}\Tfun \ar@{..>}[r] & \Efun
}
\]
whenever $\ell\ge k+8$. The unbroken arrows are obvious inclusions. It is easy to deduce a) from that diagram, with $\ell$ in place
of $k$, if $\ell$ is even. We leave that to the reader. (Recall that $T_2\Ofun$ is
rationally equivalent, over $B\ZZ/2$, to $\Ofun$
except for a possible deviation at $V=0$. Therefore $L_{[1,2]}\Ofun$
is rationally equivalent to the functor $V\mapsto \Or(V\oplus\RR^\infty)/\Or(V)$
except for a possible deviation at $V=0$.)
\newline
We can describe $L_{[1,2]}\Ofun$ as
$\hofiber[\,\kappa_2\co L_1\Ofun\to \Omega^{-1}L_2\Ofun\,]$,
where $\Omega^{-1}L_2\Ofun$ is a homogeneous functor of degree $2$, and a delooping of $L_2\Ofun$.
Similarly $L_{[1,2]}\Tfun$ can be described as $\hofiber[\,\kappa_2\co L_1\Tfun\to \Omega^{-1}L_2\Tfun\,]$
where $\Omega^{-1}L_2\Tfun$ is a homogeneous functor of degree $2$, and a delooping of $L_2\Tfun$.
By proposition~\ref{prop-ratprin} we may pretend that $L_2\Ofun$ and $L_2\Tfun$ and their deloopings
correspond to spectra with action of $\Or(2)$ whose homotopy groups are rational vector spaces. In
particular, $\Omega^{-1}L_2\Ofun$ then corresponds to $\underline{H}\QQ$ with the trivial action of $\Or(2)$.
\newline
What we are trying to do, therefore, is to produce a commutative diagram in the shape of a prism:
\[
\xymatrix@M=14pt@R=-10pt@C=-6pt{
&&& \Omega^{-1}L_2\Ofun \ar[dddlll] \ar[dddrrr] &&& \\
\\
&&& L_1\Ofun \ar[dddlll] \ar[dddrrr] \ar[uu] &&& \\
{\rule{3mm}{0mm}\Omega^{-1}L_2\Tfun}    \ar@{..}[rr] && \ar@{..}[rr] &&  \ar@{..>}[rr] &&
{\!\!\Omega^{-1}L_2\Ofun(-\oplus\RR^\ell)}  \\
\\
L_1\Tfun \ar[uu] \ar@{..>}[rrrrrr] &&&&&& L_1\Ofun(-\oplus\RR^\ell) \ar[uu]
}
\]
The three vertical arrows are of type $\kappa_2$\,. All
undotted arrows are already given and the resulting rectangle-shaped faces or parallelograms
are commutative. The front horizontal (dotted) arrow is easy to produce for any $\ell\ge 0$
(as in the proof of proposition~\ref{prop-hypc}), making the front triangle
commutative. The back horizontal (dotted) arrow is given to us for $\ell=k$, since we
are assuming a) and are aware of lemma~\ref{lem-specshift}. By assumption it makes the back triangle in the diagram commutative.
What we now have to arrange is (1) commutativity in the bottom square, and (2)
commutativity in the prism as a whole. (Making it easier using cofibrant replacements
where appropriate is allowed and recommended.)
Regarding (1), we have to connect two elements of
\[ \Rnat_\pt(L_1\Tfun,\Omega^{-1}L_2\Ofun(-\oplus\RR^\ell)) \]
by a path; these are given with $\ell=k$, but we have permission to increase $\ell$. By
lemma~\ref{lem-shifty},
a solution exists if we increase from $\ell=k$ to $\ell=k+4$. (We choose
$\ell=k+4$ rather than $\ell=k+3$ because that makes it easier to apply lemma~~\ref{lem-shifty}
a second time.) Regarding (2), we have a map from $S^1$ to
\[ \Rnat_\pt(L_1\Ofun,\Omega^{-1}L_2\Ofun(-\oplus\RR^\ell)) \]
and we have to extend that map to a disk $D^2$. (Indeed, since (1) has been taken care of,
each of the five 2-dimensional faces of the prism
determines two maps from $L_1\Ofun$ to $\Omega^{-1}L_2\Ofun(-\oplus\RR^\ell)$ and a
homotopy between those two. By joining the five homotopies, we obtain a map from a circle to
$\Rnat_\pt(L_1\Ofun,\Omega^{-1}L_2\Ofun(-\oplus\RR^\ell))$ as claimed.)
Again we have permission to increase $\ell$, and by lemma~\ref{lem-shifty}, a solution exists
if we increase from $\ell=k+4$ to $\ell=k+8$. \qed

\begin{rem} \label{rem-funnyk} We could have been more precise in stating theorem~\ref{thm-pessimist}. First,
a) for fixed even $k$ and fixed even $n\ge 4$ implies b) for the same $n$ and $k$.
Next, b) for infinitely many values of even $n$ and a specific even $k$ implies c) for
the same $k$. Finally c) for a specific even $k$ implies a) for all $n$, but with
$k+8$ instead of $k$. Hence if a) holds for infinitely many even $n$ and a fixed even $k$, then it holds
for all even $n$ with $k+8$ instead of $k$.
\end{rem}

\noindent\textit{Acknowledgments.}
This project was supported by the Engineering and Physical Sciences Research Council (UK),
Grant EP/E057128/1, and in the later stages through a Humboldt Professorship award to M. Weiss.

\end{document}